\documentclass[12pt]{amsart}
\makeindex
\usepackage{graphicx,color} 
\usepackage{amssymb}
\usepackage{hyperref} 
\usepackage[latin1]{inputenc}  

\newtheorem{propo}{Proposition}[section]
\newtheorem{theo}{Theorem}[section] 
\newtheorem{lemma}[propo]{Lemma}
\newenvironment{demo}[1]{\textit{Proof #1. }}{\hfill$\diamondsuit$}

\newtheorem{corol}[propo]{Corollary}
\newtheorem{defi}[propo]{Definition}
 
\numberwithin{equation}{section}

\let\Diam\diamondsuit

\newcommand{\surj}{\to\kern-.1em\llap{$\to$}}

\newcommand{\jrus}{\leftarrow\kern-.1em\llap{$\leftarrow$}}
\newcommand{\strictsubset}{\hbox{$\subseteq\kern-.4em\llap{${}_/$}$}}

\renewcommand{\(}{\left(}
\renewcommand{\)}{\right)}

\newcommand{\norm}[1]{\lVert#1\lVert}
\newcommand{\abs}[1]{\left\lvert#1\right\lvert}

\begin{document}
\title{{Discrete Riemann surfaces}}
\author{Christian  \textsc{Mercat}}
\email{
\href{mailto:mercat@math.univ-montp2.fr}{mercat@math.univ-montp2.fr}} 
\address{Départment de Mathématiques\\
Université Montpellier II\\ 34095 Montpellier cedex 5}
\begin{abstract}
  We detail the theory of Discrete Riemann Surfaces. It takes place on
  a cellular decomposition of a surface, together with its Poincaré
  dual, equipped with a discrete conformal structure. A lot of
  theorems of the continuous theory follow through to the discrete
  case, we will define the discrete analogs of period matrices,
  Riemann's bilinear relations, exponential of constant argument and
  series. We present the notion of criticality and its relationship
  with integrability.
\end{abstract}
\maketitle

\tableofcontents

\section{Introduction}
Riemann surfaces theory was a major achievment of  XIXth century
mathematics, setting the framework where modern complex analysis
bloomed. Nowadays, surfaces are intensively used in computer science
for numeric computations, ranging from visualization to pattern
recognition and approximation of partial differential equations. A lot
of these computations involve, at the continuous level, analytic
functions. But very few algorithms care about this specificity,
although analytic functions form a relatively small vector space among
the space of functions, problems are usually crudely discretized in a
way that doesn't take advantage of analyticity.

The theory of discrete Riemann surfaces aims at filling this gap and
setting the theoretical framework in which the notion of discrete
analyticity is set on solid grounds.


Most of the results in this paper are a straightforward application of
the continuous theory~\cite{FK,Mum} together with the results
in~\cite{M,M01,M0206041}, to which we refer for details.  We define
the discrete period matrix, which is twice as large as in the
continuous case: the periods of a holomorphic form on the graph and on
its dual are in general different, but the continuous limit theorem,
given a refining sequence of critical maps, ensures that they converge
to the same value.  The main tool is the same as in the continuous
case, the Riemann bilinear relations.

\section{Discrete Riemann surfaces} \label{sec:Definitions}
\subsection{Discrete Hodge theory}
We recall in this section basic definitions and results
from~\cite{M01} where the notion of discrete Riemann surfaces was
defined.  We are interested in discrete surfaces given by a cellular
decomposition $\diamondsuit$ of dimension two, where all faces are
\emph{quadrilaterals} (a quad-graph~\cite{Ken02,BoS, M0402097}).  
\begin{figure}[htbp]
\begin{center}\input{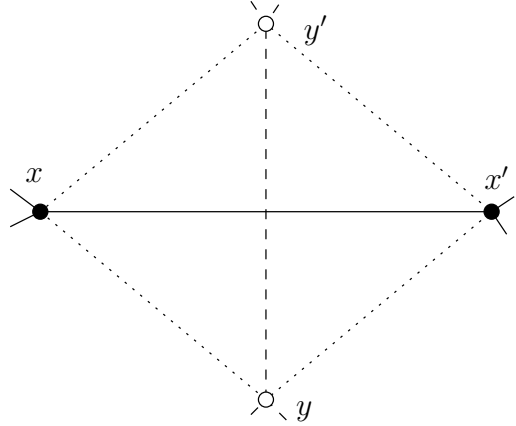}
\end{center}
\caption{The vertices and diagonals of a quadrilateral define a pair
  of dual edges.}  \label{fig:diamGG}
\end{figure}
Its vertices and diagonals define, up to
homotopy and away from the boundary, two dual\index{duality} cellular decompositions\index{cellular decomposition}
$\Gamma$ and $\Gamma^*$:
\begin{figure}[htbp]
\begin{center}\input{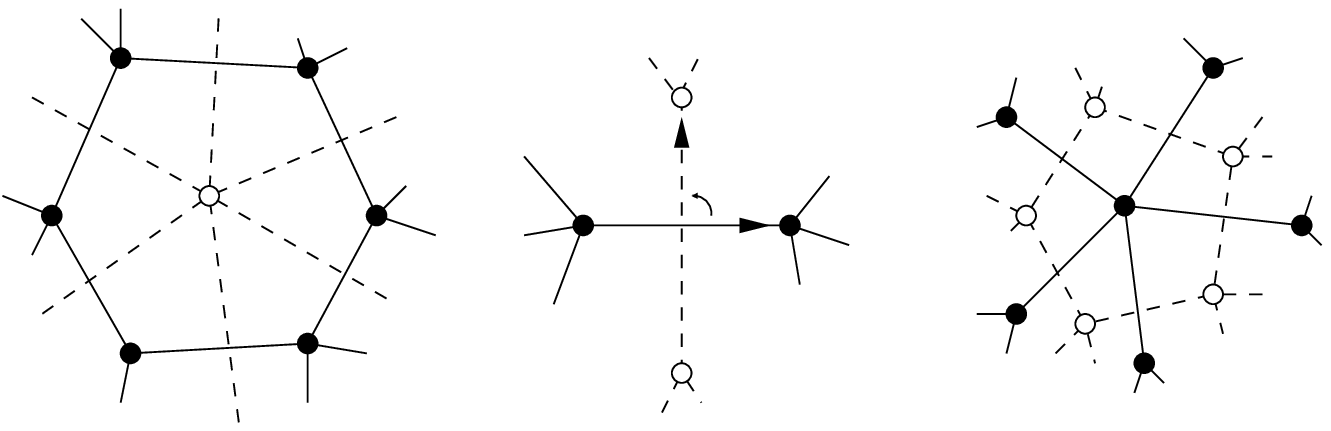}
\end{center}
\caption{Duality.}      \label{fig:duality}
\end{figure}
The edges in $\Gamma^*_1$ are dual to edges in $\Gamma_1$, faces in
$\Gamma^*_2$ are dual to vertices in $\Gamma_0$ and vice-versa.  Their
union is denoted the \emph{double} $\Lambda=\Gamma\sqcup \Gamma^*$. A
\emph{discrete conformal structure}\index{discrete conformal
  structure} on $\Lambda$ is a real positive function $\rho$ on the
unoriented edges satisfying $\rho(e^*)=1/\rho(e)$. It defines a
genuine Riemann surface structure on the discrete surface: Choose a
length $\delta$ and realize each quadrilateral by a lozenge whose
diagonals have a length ratio given by $\rho$.  Gluing them together
provides a flat riemannian metric with conic singularities at the
vertices, hence a conformal structure~\cite{Tro}.  It leads to a
straightforward discrete version of the \emph{Cauchy-Riemann
  equation}\index{discrete Cauchy-Riemmann equation}.  A function on the vertices is discrete holomorphic iff
for every quadrilateral $(x,y,x',y')\in\diamondsuit_2$,
\begin{equation}
  \label{eq:CR}
  f(y')-f(y)=i\,\rho(x,x')\left( f(x')-f(x)\right).
\end{equation}
\begin{figure}[htbp]
\begin{center}\input{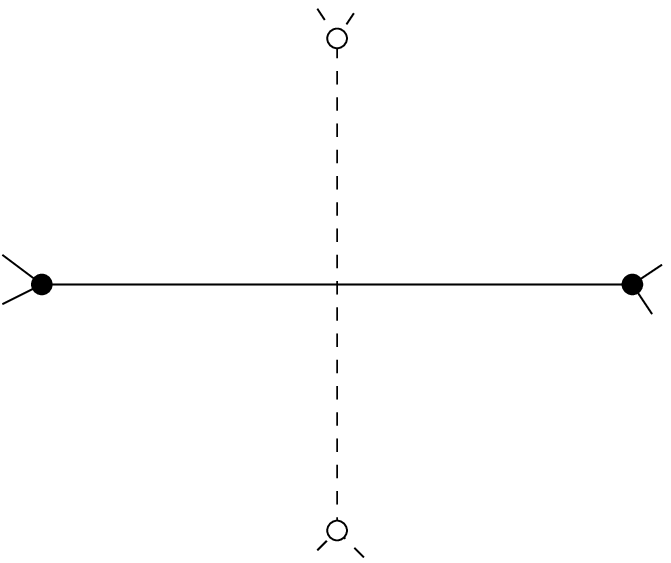}
\end{center}
\caption{The discrete Cauchy-Riemann equation.}         \label{fig:CR}
\end{figure}

We recall elements of de-Rham cohomology, doubled in our context: The complex
of \emph{chains} $C(\Lambda)=C_{0}(\Lambda)\oplus C_{1}(\Lambda)\oplus
C_{2}(\Lambda)$ is the vector space span by vertices, edges and faces.  It is
equipped with a \emph{boundary} operator $\partial:C_{k}(\Lambda)\to
C_{k-1}(\Lambda)$, null on vertices and fulfilling $\partial^{2}=0$.  The
kernel $\text{ker~}\partial=:Z_{\bullet}(\Lambda)$ of the boundary operator
are the closed chains or \emph{cycles}.  Its image are the \emph{exact}
chains.  It provides the dual spaces of forms, called \emph{cochains},
$C^{k}(\Lambda):=\text{Hom}(C_{k}(\Lambda),\mathbb{C})$ with a
\emph{coboundary} $d:C^k(\Lambda)\to C^{k+1}(\Lambda)$ defined by  Stokes
formula: $$\int\limits_{(x,x')} df:= f\left(\partial(x,x')\right)=f(x')-f(x),
\qquad \iint\limits_F d\alpha:=\oint\limits_{\partial F}\alpha.
$$
A \emph{cocycle} is a closed cochain and we note $\alpha\in 
Z^{k}(\Lambda)$. 

These spaces are equipped with the canonical scalar product, weigthed
according to $\rho$ on edges and averaged on the graph and its dual:
$$(\alpha,\beta):=\tfrac12 \sum_{e\in\Lambda_1}
\rho(e)\left(\int_e \alpha\right)\left(\int_e \beta\right).$$

Duality of complexes allows us to define a \emph{Hodge operator}\index{discrete Hodge star}
$*$ on forms by
\begin{eqnarray}
    *:C^{k}(\Lambda) & \to & C^{2-k}(\Lambda)
    \notag  \\
   \phantom{*:}C^{0}(\Lambda)\ni f & \mapsto & * f : \iint_{F}*f := f(F^{*}),
    \notag  \\
     \phantom{*:}C^{1}(\Lambda)\ni\alpha & \mapsto & * \alpha :
     \int_{e}*\alpha := -\rho(e^{*})\int_{e^{*}}\alpha,
    \label{eq:*Def}  \\
     \phantom{*:}C^{2}(\Lambda)\ni \omega & \mapsto & * \omega : 
     (*\omega)(x) := \iint_{x^{*}}\omega.
    \notag  
\end{eqnarray}

It fulfills $*^{2}=(-\text{Id}_{C^k})^{k}$. The endomorphism
$\Delta:=-d*d*-*d*d$ is the usual discrete \emph{Laplacian}\index{discrete Laplacian}: Its
formula on a function at a vertex $x\in\Gamma_0$ with neighbours
$x_1,\ldots,x_V\in\Gamma_0$ is the usual weighted averaged difference:
$$\left(\Delta(f)\right)(x)=
\sum_{k=1}^V\rho(x,x_k)\left(f(x)-f(x_k)\right).$$
 The space of \textit{harmonic forms} is defined as its kernel .

 The Hodge star and the laplacian are real operators. Since
 $*^2=-\text{Id}$ on functions, it is natural to consider them on
 complexified cochains. The discrete holomorphic forms are special
 complex harmonic forms: a $1$-form
\begin{equation}
\alpha\in C^{1}(\Lambda) \text{~~is \emph{holomorphic} iff~~}
d\alpha=0 \text{~~and~~} *\alpha = -i\alpha, 
\label{eq:holoDef}
\end{equation}
that is to say if it is closed and of type $(1,0)$. Let $d'$, resp.
$d''$ the compositions of the exterior derivative with the projection
on the space of $(1,0)$, resp. $(0,1)$-forms, Eq.~\eqref{eq:holoDef}
is equivalent to $d'\alpha=0$.  We will note
$\alpha\in\Omega^{1}(\Lambda)$.  A function
$f:\,\Lambda_{0}\to\mathbb{C}$ is \emph{holomorphic} iff $df$ is
holomorphic, which is equivalent to \eqref{eq:CR} and we note
$f\in\Omega^{0}(\Lambda)$.

In the compact case, $-*\,d\,*$ is the adjoint $d^*$ of the coboundary
operator $d$ and the Hodge theorem orthogonally decomposes forms into
exact, coexact and harmonic,
$$C^k(\Lambda) = \text{Im } d\oplus^\perp
 \text{Im } d^*\oplus^\perp  \text{Ker } \Delta,$$
 harmonic forms are the closed and
co-closed ones, and harmonic $1$-form are the orthogonal sum of
holomorphic and anti-holomorphic ones:
$$ \text{Ker }\Delta= \text{Ker }d\cap\text{Ker }d^*
=\text{Ker }d '\oplus^\perp\text{Ker }d''.$$

\subsection{Wedge product}
We construct a wedge product on $\Diam$ such that 
\begin{itemize}
\item the canonical weighted hermitian scalar product reads as expected
$$(\alpha,\beta)=\iint \alpha\wedge\,*\bar\beta,$$
\item and the coboundary operator $d_{\Diam}$ on $\Diam$, is a derivation
  for this product $\wedge:C^k(\Diam)\times C^l(\Diam)\to
  C^{k+l}(\Diam)$.
\end{itemize}
It is defined
by the following formulae, for $f,g\in C^0(\Diam)$, $\alpha,\beta\in
C^1(\Diam)$ and $\omega\in C^2(\Diam)$:
\begin{align}
  (f\cdot g)(x):=&f(x)\cdot g(x)\qquad \mathrm{ ~for~} x\in\Diam_0,
\label{eq:wedge00}\\ 
  \int\limits_{(x,y)}f\cdot\alpha:=& \frac{f(x)+f(y)}2
  \int\limits_{(x,y)}\alpha\qquad \mathrm{~for~} (x,y)\in\Diam_1,
\label{eq:wedge01}\\ 
  \iint\limits_{\hidewidth{(x_1,x_2,x_3,x_4)}\hidewidth}\alpha\wedge\beta
:=&\tfrac{1}{4}\sum_{k=1}^4
\int\limits_{{(x_{k-1},x_k)}}\alpha\;
\int\limits_{\hidewidth{(x_k,x_{k+1})}\hidewidth}\beta-
\int\limits_{{(x_{k+1},x_k)}}\alpha\;
\int\limits_{\hidewidth{(x_k,x_{k-1})}\hidewidth}\beta,
\label{eq:wedge11}\\  
\iint\limits_{\hidewidth{(x_1,x_2,x_3,x_4)}\hidewidth}
  f\cdot\omega:=&\frac{\scriptstyle f(x_1)+f(x_2)+f(x_3)+f(x_4)}{4}
  \iint\limits_{\hidewidth{(x_1,x_2,x_3,x_4)}\hidewidth}\omega
\label{eq:wedge02}\\ 
  &\qquad \mathrm{~for~} (x_1,x_2,x_3,x_4)\in\Diam_2.\notag
\end{align}

\begin{figure}[htbp]
\begin{center}\input{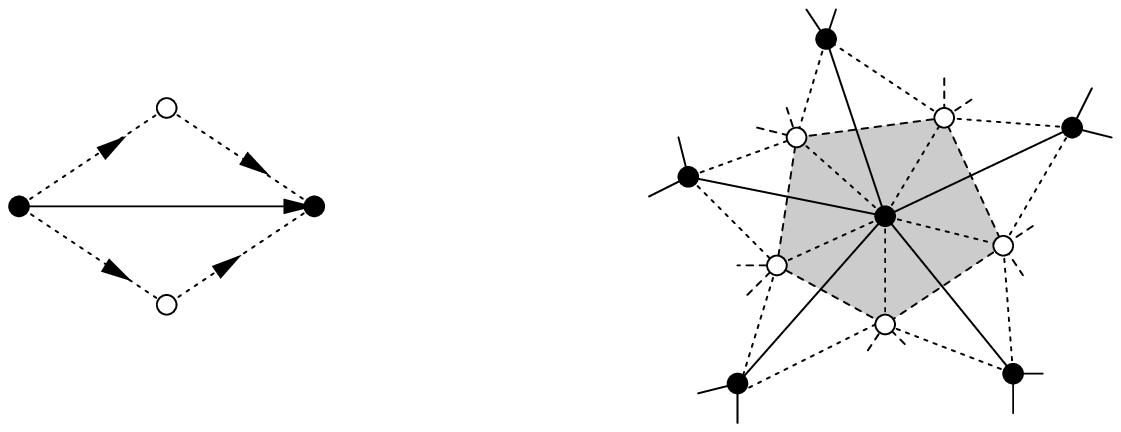}
\end{center}
\caption{Notations.}    \label{fig:avera}
\end{figure}

A form on $\Diam$ can be \emph{averaged} into a form on $\Lambda$:
This map $A$ from 
$C^\bullet(\Diam)$
to $C^\bullet(\Lambda)$ is the identity for functions and
defined by the following formulae for $1$ and $2$-forms:
\begin{align}
 \int\limits_{{(x,x')}}A(\alpha_\Diam):= \tfrac12 \left(
   \int\limits_{{(x,y)}}+\int\limits_{\hidewidth{(y,x')}\hidewidth}
   +\int\limits_{\hidewidth{(x,y')}\hidewidth}+\int\limits_{{(y',x')}}
 \right) \alpha_\Diam, \label{def:avera1}\\ 
 \iint\limits_{\hidewidth{x^*}\hidewidth}A(\omega_\Diam):=
\tfrac12\sum_{k=1}^d\;\iint\limits_{{(x_k,y_k,x,y_{k-1})}}\omega_\Diam,
\label{def:avera2}
\end{align}
where notations are made clear in Fig.~\ref{fig:avera}.  The map $A$
is neither injective nor surjective in the non simply-connected case,
so we can neither define a Hodge star on $\Diam$ nor a wedge product
on $\Lambda$.  Its kernel on $1$-forms is $\text{\rm Ker
}(A)=\text{\rm Vect }(d_\Diam\varepsilon)$, where $\varepsilon$ is the
biconstant, $+1$ on $\Gamma$ and $-1$ on $\Gamma^*$.  But $d_\Lambda
A=A\,d_\Diam$ so it carries cocycles on $\Diam$ to cocycles on
$\Lambda$. Its image are these cocycles of $\Lambda$ verifying that
their holonomies along cycles of $\Lambda$ only depend on their
homology on the combinatorial surface. Given a $1$-cocycle $\mu\in
Z^{1}(\Lambda)$ with such a property, a corresponding $1$-cocycle
$\nu\in Z^{1}(\Diam)$ is built in the following way: Choose an edge
$(x_{0},y_{0})\in\Diam_{1}$; for an edge $(x,y)\in\Diam_{1}$ with $x$
and $x_{0}$ on the same leaf of $\Lambda$, choose two paths
$\lambda_{x,x_{0}}$ and $\lambda_{y_{0},y}$ on the double graph
$\Lambda$, from $x$ to $x_{0}$ and $y_{0}$ to $y$ respectively, and
define \begin{equation}
\int_{(x,y)}\nu:=\int_{\lambda_{x,x_{0}}}\mu+\int_{\lambda_{y_{0},y}}\mu
-\oint_{[\gamma]}\mu \label{eq:holonomyProp} \end{equation} where
$[\gamma]=[\lambda_{x,x_{0}}+(x_{0},y_{0})+\lambda_{y_{0},y}+(y,x)]$
is the class of the full cycle in the homology of the
surface. Changing the base points change $\mu$ by a multiple of
$d_\diamondsuit\varepsilon$.

It follows in the compact case that the dimensions of the harmonic
forms on $\Diam$ (the kernel of $\Delta A$) modulo $d\varepsilon$, as
well as the harmonic forms on $\Lambda$ with same holonomies on the
graph and on its dual, are twice the genus of the surface, as
expected.  Unfortunately, the space $\text{Im~}
A=\mathcal{H}^\perp\oplus\text{Im }d$ is not stable by the Hodge star
$*$.  We could nevertheless define holomorphic $1$-forms on $\Diam$ but
their dimension would be much smaller than in the continuous, namely the
genus of the surface.  Criticality provides conditions which ensure
that the space $*\text{Im }A$ is ``close'' to $\text{Im }A$.

We 
construct an
{\em heterogeneous}\index{discrete wedge product} wedge product for $1$-forms: with $\alpha,\beta\in
C^1(\Lambda)$, define $\alpha\wedge\beta\in C^1(\Diam)$ by
\begin{equation}
    \iint\limits_{\hidewidth(x,y,x',y')\hidewidth}\alpha\wedge\beta:=
\tfrac12\left(\;\;
    \int\limits_{\hidewidth(x,x')\hidewidth}\alpha
    \int\limits_{\hidewidth(y,y')\hidewidth}\beta+
    \int\limits_{\hidewidth(y,y')\hidewidth}\alpha
    \int\limits_{\hidewidth(x',x)\hidewidth}\beta
\right).
    \label{eq:averaging}
\end{equation}

It verifies $A(\alpha_\Diam)\wedge A(\beta_\Diam)=\alpha_\Diam\wedge
\beta_\Diam$, the first wedge product being between $1$-forms on
$\Lambda$ and the second between forms on $\Diam$.  The usual scalar
product on compactly supported forms on $\Lambda$ reads
as expected:
 \begin{equation}
     (\alpha,\beta)=\tfrac12\sum_{e\in\Lambda_1}\rho(e)\left(
     \int_e\alpha\right) \left(
     \int_e\bar\beta\right)=\iint\limits_{\Diam_2}\alpha\wedge*\bar\beta
     \label{eq:scalarProd}
 \end{equation}

 \subsection{Energies}
 The $L^2$ norm of the $1$-form $df$, called the Dirichlet energy of
 the function $f$, is the mean of the usual Dirichlet energies on
 each independant graph:
\begin{align}
  E_D(f):=\tfrac12 \lVert df\rVert^2 &=
  \tfrac12 \left(df,\,df\right)=\tfrac14
\sum_{(x,x')\in\Lambda_1}\rho(x,x')
  \left\lvert f(x') - f(x) \right\rvert^2\label{eq:norm}
  \\
  &=\frac{ E_D(f|_\Gamma)+ E_D(f|_{\Gamma^*})}2
.\notag\end{align}
Harmonic maps minimize this energy among functions fulfilling certain
boundary conditions.

The conformal energy of a map measures its conformality defect, it is
null on holomorphic functions:
\begin{equation}
  E_C(f) := \tfrac14\lVert df -i *  df\rVert^2.
\label{eq:EC}
\end{equation}
It is related to the Dirichlet energy through the same formula as in
the continuous:
\begin{align}
E_C(f) &= \tfrac14\left( df  -i * df,\, df  -i * df\right)
\notag\\
&= \tfrac14\lVert df \rVert^2+\tfrac14\lVert-i * df \rVert^2
+\, \tfrac12\text{Re}(df,\, -i * df)\notag\\
&=
  \tfrac12\lVert df \rVert^2 + \, \tfrac12\text{Im} \iint_{\lozenge_2} df\wedge\overline{df}
\notag\\
&=
E_D(f) - \mathcal{A}(f)
\label{eq:ECEDA}
\end{align}
where the area of the image of the application $f$ in the complex
plane has the same formula 
\begin{equation}
\mathcal{A}(f) = \tfrac i2
\iint_{\lozenge_2} df\wedge\overline{df}
\label{eq:A}
\end{equation}
 as in the continuous case
since, for a face $(x,y,x',y')\in\lozenge_2$, the algebraic area of the
oriented quadrilateral $\Bigl(f(x),f(x'),f(y),f(y')\Bigr)$ is given by
\begin{align*}
  \smash{\iint\limits_{(x,y,x',y')}} df\wedge\overline{df}&=
i\, \text{Im}\left(
(f(x')-f(x))\overline{(f(y')-f(y))}
\right)\\
&=-2 i\mathcal{A}\Bigl(f(x),f(x'),f(y),f(y')\Bigr).
\end{align*}

\section{Period matrix} \label{sec:PeriodMatrix}
We use the convention of Farkas and Kra~\cite{FK}, chapter III, to
which we refer for details.  Consider $(\Diam,\rho)$ a discrete
compact Riemann surface.

\subsection{Intersection number, on $\Lambda$ and on $\diamondsuit$} \label{sec:Intersection}
For a given simple (real) cycle $C\in Z_{1}(\Lambda)$, we construct a
harmonic $1$-form $\eta_{C}$ such that $\oint_{A}\eta_{C}$ counts the
algebraic number of times $A$ contains an edge dual to an edge of $C$:
It is the solution of a Neumann problem on the surface cut open along
$C$ (see~\cite{M} for details).  It follows from standard homology
technique that $\eta_{C}$ depends only on the homology class of $C$
(all the cycles which differ from $C$ by an exact cycle $\partial A$)
and can be extended linearly to all cycles as $\eta_\bullet:
H_1(\Lambda)\to C^1(\Lambda)$; it fulfills, for a closed form
$\theta$,
\begin{equation}
  \label{eq:dualityEta}
  \oint_C \theta=\iint_\diamondsuit \eta_C\wedge\theta,
\end{equation}
and a basis of the homology provides a dual basis of harmonic forms on
$\Lambda$.  Beware that if the cycle $C\in Z_{1}(\Gamma)$ is purely on
$\Gamma$, then this form ${\eta_{C}}_{|_\Gamma}=0$ is null on $\Gamma$.
 
\begin{figure}[htbp]
\begin{center}\input{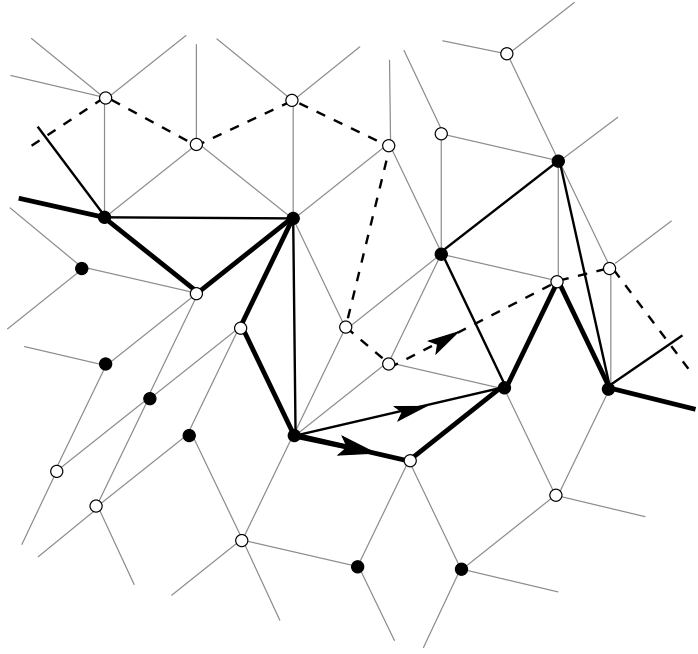}
\end{center}
\caption{A path $C$ on $\diamondsuit$ defines a pair of paths $C_{\Gamma}$ and  $C_{\Gamma^*}$ on its left.}
\label{fig:diamGamGams}
\end{figure}

The \emph{intersection number}\index{intersection number} between two cycles $A,B\in 
Z_{1}(\Lambda)$ is defined as
\begin{equation}
    A\cdot B:=\iint_{\Diam}\eta_{A}\wedge\eta_{B}.
    \label{eq:interDef}
\end{equation}
It is obviously linear and antisymmetric, it is an integer number for integer
cycles. Let's stress again that the intersection of a cycle on $\Gamma$ with
another cycle on $\Gamma$ is always null.  A cycle $C\in Z_{1}(\Diam)$
defines a pair of cycles on each graph $C_{\Gamma}\in Z_{1}(\Gamma)$,
$C_{\Gamma^{*}}\in Z_{1}(\Gamma^{*})$ which are homologous to $C$ on the
surface, composed of portions of the boundary of the faces on $\Lambda$ dual
to the vertices of $C$.  They are uniquely defined if we require that they
lie ``to the left'' of $C$ as shown in Fig.\ref{fig:diamGamGams}.
By the procedure~\eqref{eq:holonomyProp} applied to
$\eta_{C_{\Gamma}}+\eta_{C_{\Gamma^{*}}}$, we construct a $1$-cocycle
$\eta_{C}\in Z^{1}(\Diam)$ unique up to $d\varepsilon$, and since
$\forall\theta,\; d\varepsilon\wedge\theta=0$, Eq.~\eqref{eq:interDef}
defines an intersection number on $Z_{1}(\Diam)$.  Unlike the intersection
number on $\Lambda$, this one has all the usual expected properties. In
particular Eq.~\eqref{eq:interDef} holds for $A,B\in Z_{1}(\diamondsuit)$.

\subsection{Canonical dissection, fundamental polygon} \label{sec:CanonicalDiss}
The complex $\Diam$ being connected, consider a maximal tree
$T\subset\Diam_{1}$, that is to say $T$ is a
$\mathbf{Z}_2$-homologically trivial chain and every edge added to $T$
forms a cycle.  A \emph{canonical dissection}\index{canonical
  dissection} or cut-system $\aleph$ of the genus $g$ discrete Riemann
surface $\Diam$ is given by a set of oriented edges $(e_{k})_{1\leq
  k\leq 2g}$ such that the cycles $\aleph\subset(T\cup e_{k})$ form a
basis of the homology group $H_{1}(\Diam)$ verifying, for $1\leq
k,\ell\leq g$
\begin{equation}
    \aleph_{k}\cdot\aleph_{\ell}=0,   \quad
    \aleph_{k+g}\cdot\aleph_{\ell+g}=0,   \quad
    \aleph_{k}\cdot\aleph_{\ell+g}=\delta_{k,\ell}.
 \label{eq:alephDiam} 
\end{equation}
They actually form a basis of the fundamental group
$\pi_1(\diamondsuit)$ and the defining relation among them is (noted
multiplicatively)
\begin{equation}
  \label{eq:pi1}
  \prod_{k=1}^g\aleph_{k}\aleph_{k+g}\aleph_{k}^{-1}\aleph_{k+g}^{-1}=1.
\end{equation}
The construction of such a basis is standard and we won't repeat the
procedure. What is less standard is the interpretation of
Eq.~\eqref{eq:pi1} in terms of the boundary of a fundamental domain,
discretization introduces some subtleties (that can safely be skept in
first instance). We end up with the familiar $2g\times 2g$
intersection numbers matrix on $\diamondsuit$.

\begin{figure}[htbp]
\begin{center}\input{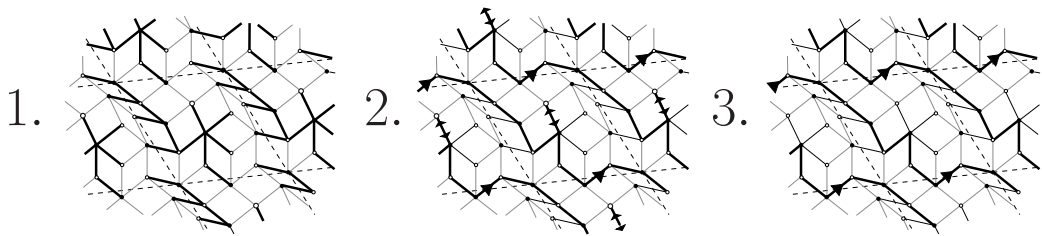}
\end{center}
\caption{1. A maximal rooted tree in a quadrilateral decomposition of the
  torus. 2. An additional edge defines a rooted cycle $\aleph_1$, pruned of
  its dangling trees. 3. Its un-rooted version, the simple loop $\aleph_1^-$.}
\label{fig:prune}
\end{figure}

Considering $T\cup e_{k}$ as a rooted graph, we can prune it of all its
pending branches, leaving a simple closed loop $\aleph_k^-$, attached to the
origin $O$ by a simple path $\lambda_k$ (see Fig.~\ref{fig:prune}), yielding
the cycle $\aleph_k$. These three cycles are deformation retract of one
another, $\aleph_k^-\subset \aleph_k\subset T\cup e_{k}$ hence are equal in
homology.

In the continuous case, a basis of the homology can be realized by $2g$
simple arcs, transverse to one another and meeting only at the base point. It
defines an isometric model of the surface as a fundamental domain
homeomorphic to a disc and bordered by $4g$ arcs to identify pairwise.  In
the discrete case, by definition, the set $\diamondsuit\setminus\aleph$ of
the cellular complex minus the edges taking part into the cycles basis is
homeomorphic to a disc hence the surface is realized as a polygonal
\emph{fundamental domain} $\mathcal{M}$ whose boundary edges are identified
pairwise.

But it is sometimes impossible to choose a basis of the homology
verifying~\eqref{eq:alephDiam} by simple discrete cycles which are
transverse to one another. For instance, if the path $\lambda_k$ is
not empty, the cycle $\aleph_k$ is not even simple. Moreover, some
edges may belong to several cycles. In this case, the edges on the
boundary of this fundamental polygon can not be assigned a unique
element of the basis or its inverse, and therefore can not be grouped
into only $4g$ continuous paths to identify pairwise but more than
$4g$.

\begin{figure}[htbp]
\begin{center}\input{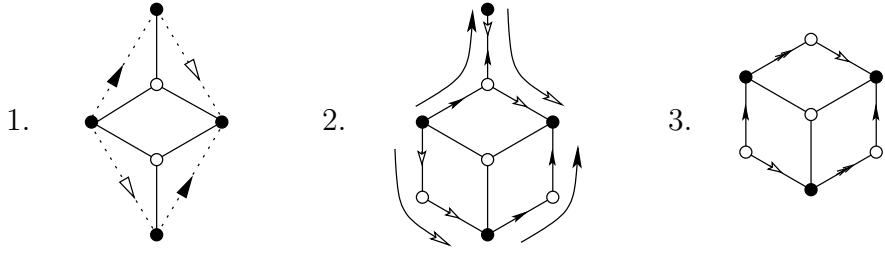}
\end{center}
\caption{Three different fundamental polygons of a decomposition of the torus
  ($g=1$) by three quadrilaterals: 1. The standard fundamental domain where
  the $4g$ paths are not adapted to $\diamondsuit$. 2. $\mathcal{M}_+$ is
  composed of edges of $\diamondsuit$ composing $4g$ arcs (which may have
  portions in common) to identify pairwise, each edge corresponds to an
  element of the basis $\aleph$ or its inverse, except for edges of
  ``dangling trees'' which are associated with two such elements.  3.
  $\mathcal{M}$ is composed of edges of $\diamondsuit$ composing more than
  $4g$ arcs to identify pairwise, there is no correspondence with a basis of
  cycles.}
\label{fig:FundPoly}
\end{figure}

In fact, the information contained into the basis $\aleph$ is more than
simply this polygon, the set of edges composing the concatenated cycle
\begin{equation}
  \label{eq:FundPoly1}
(\aleph_1,\aleph_{g+1},\aleph_1^{-1},\aleph_{g+1}^{-1},\aleph_2,\ldots,
\aleph_{g}^{-1},\aleph_{2g}^{-1})
\end{equation}
encodes a cellular complex $\mathcal{M}_+$ which is \emph{not} a
combinatorial surface and consists of the fundamental polygon $\mathcal{M}$
plus some \emph{dangling trees}, corresponding to the edges which belong to
more than one cycle or participate more than once in a cycle (the paths
$\lambda_k$), as exemplified in Fig.\ref{fig:FundPoly}. By construction, the
edge $e_k$ belongs to the cycle $\aleph_k$ only, hence these trees are in
fact without branches, simple paths whose only leaf is the base point $O$. To
retrieve the surface, the edges of this structure $\mathcal{M}_+$ are
identified group-wise, an edge participating $k$ times in cycles will have
$[k/2]+2$ representatives to identify together, two on the fundamental
polygon and the rest as edges of dangling trees.

Eliminating repetition, that is to say looking at~\eqref{eq:FundPoly1} not as
a \emph{sequence} of edges but as a simplified \emph{cycle} (or a simplified
word in edges), thins $\mathcal{M}_+$ into $\mathcal{M}$, pruning away the
dangling paths.  The fundamental polygon boundary loses its structure as $4g$
arcs to be identified pairwise, in general a basis cycle will be disconnected
around the fundamental domain and a given edge can not be assigned to a
particular cycle.  This peculiarity gives a more complex yet well defined
meaning to the contour integral formula for a $1$-form $\theta$ defined on
the boundary edges of $\mathcal{M}_+$,
\begin{equation}
    \oint_{\partial\mathcal{M}}\theta=\sum_{k=1}^{2g}
    \oint_{\aleph_{k}}\theta+\oint_{\aleph_{k}^{-1}}\theta.
    \label{eq:contourInt}
\end{equation}
This basis gives rise to cycles $\aleph^{\Gamma}$ and $\aleph^{\Gamma^{*}}$
whose homology classes form a basis of the group for each respective graph,
that we compose into $\aleph^{\Lambda}$ defined by
\begin{align}
    \aleph^{\Lambda}_{k} & =  \aleph^{\Gamma}_{k},
      &
    \aleph^{\Lambda}_{k+g} & =  \aleph^{\Gamma^{*}}_{k},
    \label{eq:alephLambda}  \\
    \aleph^{\Lambda}_{k+2g} & =  \aleph^{\Gamma^{*}}_{k+g},
     &
    \aleph^{\Lambda}_{k+3g} & =  \aleph^{\Gamma}_{k+g},
    \notag
\end{align}
for $ 1\leq k\leq g $ so that while the intersection numbers matrix on
$\diamondsuit$ is given by the $2g\times 2g$ matrix
\begin{equation}
    (\aleph^{}_{k}\cdot\aleph^{}_{\ell})_{k,\ell}=
    \begin{pmatrix}
        0 & I  \\
        -I & 0  
    \end{pmatrix},
    \label{eq:alephIntersecDiam}
\end{equation} 
the intersection numbers matrix on $\Lambda$ is the $4g\times 4g$ matrix with
the same structure
\begin{equation}
    (\aleph^{\Lambda}_{k}\cdot\aleph^{\Lambda}_{\ell})_{k,\ell}=
    \begin{matrix}
      \begin{matrix}
       \phantom{-}\Gamma^{ } &\phantom{-}\Gamma^* &
\!\Gamma^* & \!\Gamma^{ }
      \end{matrix}\\
      \begin{pmatrix}
0&0&
      I & 0  \\
0&0&
      0 & I
\\
     -I & 0   &0&0\\
      0 &  -I   &0&0
  \end{pmatrix}&
  \begin{matrix}
    \Gamma^{\phantom{*}}\\
    \Gamma^*\\
    \Gamma^*\\
    \Gamma^{\phantom{*}}
  \end{matrix}
  \end{matrix}
.
    \label{eq:alephIntersec}
\end{equation}

\subsection{Bilinear relations} \label{sec:Bilin}\index{bilinear relations}
 
\begin{propo} \label{prop:bilin}
     Given a canonical dissection $\aleph$, for two closed
    forms $\theta,\theta'\in Z^{1}(\Diam)$,
    \begin{equation}
        \iint_{\Diam}\theta\wedge\theta' = \sum_{j=1}^{g}\left(
        \oint_{\aleph_{j}}\theta\oint_{\aleph_{j+g}}\!\!\!\!\theta'-
        \oint_{\aleph_{j+g}}\!\!\!\!\theta\,\oint_{\aleph_{j}}\theta' \right);
        \label{eq:iintThetaDiam}
    \end{equation}
    for two closed
   forms $\theta,\theta'\in Z^{1}(\Lambda)$,
   \begin{equation}
       \iint_{\Diam}\theta\wedge\theta' = \sum_{j=1}^{2g}\left(
       \oint_{\aleph_{j}^{\Lambda}}\theta
       \oint_{\aleph_{j+2g}^{\Lambda}}\!\!\!\!\theta'-
       \oint_{\aleph_{j+2g}^{\Lambda}}\!\!\!\!\theta\,
       \oint_{\aleph_{j}^{\Lambda}}\theta' \right).
       \label{eq:iintTheta}
   \end{equation}
\end{propo}
\begin{demo}{\ref{prop:bilin}}
  Each side is bilinear and depends only on the cohomology classes of the
  forms.  Decompose the forms onto the cohomology basis $(\alpha_{k})$.  On
  $\Lambda$, use Eq~\eqref{eq:alephIntersecAlpha} for the LHS and the duality
  property Eq.~\eqref{eq:dualBasis} for the RHS. On $\diamondsuit$, use their
  counterparts.
\end{demo}

Notice that for a harmonic form $\theta\in\mathcal{H}^{1}(\Lambda)$,
the form $*\theta$ is closed as well, therefore its norm is given by
\begin{equation}
    \theta\in\mathcal{H}^{1}(\Lambda)\implies 
    \norm{\theta}^{2}= \sum_{j=1}^{2g}\left(
    \oint_{\aleph_{j}}\theta\oint_{\aleph_{j+2g}}\!\!\!\!*\bar\theta-
    \oint_{\aleph_{j+2g}}\!\!\!\!\theta\,\oint_{\aleph_{j}}*\bar\theta \right).
    \label{eq:iintThetaHarmo}
\end{equation}

\subsection{Basis of harmonic forms, basis of holomorphic forms} \label{sec:Basis}
We define $\alpha^{\Lambda}$, the basis of real harmonic $1$-forms,
dual to the homology basis $\aleph^{\Lambda}$, as described in
Sec.~\ref{sec:Intersection},
\begin{eqnarray}
\alpha^{\Lambda}_{k}&:=&\eta_{\aleph^{\Lambda}_{k+2g}} \qquad 
\text{ ~and}\notag\\
\alpha^{\Lambda}_{k+2g}&:=&-\eta_{\aleph^{\Lambda}_{k}} \qquad 
\text{ ~for~ } 1\leq k\leq 2g
\label{eq:dualBasisDef}
\end{eqnarray}
which verify
\begin{eqnarray}
\oint_{\aleph^{\Lambda}_{k}}\alpha_{\ell}&=&\delta_{k,\ell},\notag\\
\oint_{\aleph^{\Lambda}_{k+2g}}\alpha_{\ell+2g}&=&\delta_{k,\ell},
\label{eq:dualBasis}
\end{eqnarray}
and dually, the intersection matrix elements are given by
\begin{equation}
    \aleph^{\Lambda}_{k}\cdot\aleph^{\Lambda}_{\ell}= 
    \iint_{\Diam}\alpha^{\Lambda}_{k}\wedge\alpha^{\Lambda}_{\ell}=
    (\alpha^{\Lambda}_{k},-*\alpha^{\Lambda}_{\ell}).  
    \label{eq:alephIntersecAlpha}
\end{equation}
On $\Diam$, the elements $\alpha^{\Diam}_{k}:=\eta_{\aleph_{k+g}}$ and
$\alpha^{\Diam}_{k+g}:=-\eta_{\aleph_{k}}$ for $1\leq k\leq g$, defined up to
$d\varepsilon$, verify
$A(\alpha^{\Diam}_{k})=\alpha^{\Lambda}_{k}+\alpha^{\Lambda}_{k+g}$,
$A(\alpha^{\Diam}_{k+g})=\alpha^{\Lambda}_{k+2g}+\alpha^{\Lambda}_{k+3g}$ and
form a basis of the cohomology on $\Diam$ dual to $\aleph$ as well,
\begin{eqnarray}
\alpha^{\diamondsuit}_{k}&:=&\eta_{\aleph^{\diamondsuit}_{k+g}} \qquad 
\text{ ~and}\notag\\
\alpha^{\diamondsuit}_{k+g}&:=&-\eta_{\aleph^{\diamondsuit}_{k}} \qquad 
\text{ ~for~ } 1\leq k\leq g,
\label{eq:dualBasisDefDiam}
\end{eqnarray}
they fulfill the first identity in~Eq.\eqref{eq:alephIntersecAlpha} but the
second is meaningless in general since $*$ can not be defined on
$\diamondsuit$. We will drop the mention $\Lambda$ when no confusion is
possible.

\begin{propo} \label{prop:positive}
  The matrix of inner products on $\Lambda$,
 \begin{equation}
     (\alpha_{k},\alpha_{\ell})_{k,\ell}
     =  \iint_{\Diam}\alpha_{k}\wedge*\bar\alpha_{\ell}
     = \begin{cases}
        +\oint_{\aleph_{k+2g}} *\alpha_\ell,& 1\leq k\leq 2g,\\
        -\oint_{\aleph_{k-2g}} *\alpha_\ell,& 2g< k\leq 4g.
       \end{cases}
   =:
     \begin{pmatrix}
       A & D  \\
       B & C
     \end{pmatrix}
     \label{eq:GammaMatrix}
 \end{equation}  
  is a real symmetric positive definite matrix.
\end{propo}
\begin{demo}{\ref{prop:positive}}
  It is real because the forms are real, and symmetric because the scalar
  product~\eqref{eq:scalarProd} is skew symmetric. Definition
  Eq.~\eqref{eq:dualBasisDef} and Eq.~\eqref{eq:dualityEta} lead to the
  integral formulae.  Positivity follows from the bilinear
  relation~Eq.~\eqref{eq:iintTheta}: for
  $\theta=\sum_{k=1}^{4g}\xi_k\,\alpha_k$, with
  $\xi_k\in\mathbf{C},\;\sum_{k=1}^{4g}|\xi_k|^2>0$,
\begin{eqnarray}
  \norm{\theta}^2&=&\sum_{j=1}^{2g}\left[
\int_{\aleph_j}\theta\int_{\aleph_{2g+j}}*\bar\theta-\int_{\aleph_{2g+j}}\theta\int_{\aleph_{2j}}*\bar\theta\right]\notag\\
&=&\sum_{k,\ell=1}^{4g}\xi_k\,\bar\xi_\ell\sum_{j=1}^{2g}\left[
\int_{\aleph_j}\alpha_k\int_{\aleph_{2g+j}}*\alpha_\ell
-\int_{\aleph_{2g+j}}\alpha_k\int_{\aleph_{2j}}*\alpha_\ell\right]\notag\\
&=&\sum_{k,\ell=1}^{4g}\xi_k\,\bar\xi_\ell\,(\alpha_{k},\alpha_{\ell})>0.
  \label{eq:positive}
\end{eqnarray}
\end{demo}

The form $\alpha_k$ is supported by only one of the two graphs $\Gamma$ or
$\Gamma^{*}$, the form $* \alpha_k$ is supported by the other one, and the
wedge product $\theta_\Gamma\wedge\theta'_\Gamma=0$ is null for two $1$-forms
supported by the same graph. Therefore the matrices $A$ and $C$ are $g\times
g$-block diagonal and $B$ is anti-diagonal. 
\begin{equation}
  \label{eq:blockDiag}
  A=
  \begin{pmatrix}
    A_\Gamma&0\\ 0&A_{\Gamma^*}
  \end{pmatrix},\qquad
  B=
  \begin{pmatrix}
   0&B_{\Gamma^*,\Gamma}\\ B_{\Gamma,\Gamma^*}&0
  \end{pmatrix},\qquad
  C=
  \begin{pmatrix}
   C_{\Gamma^*}&0\\ 0&C_{\Gamma}
  \end{pmatrix}.
\end{equation}
The matrices of intersection numbers~\eqref{eq:alephIntersec} and of inner
products differ only by the Hodge star $*$. Because $*$ preserves harmonic
forms and the inner product, we get its matrix representation in the basis
$\alpha$,
\begin{equation}
  \label{eq:*Mat}
  *=
  \begin{pmatrix}
    -D&A\\
    -C&B
  \end{pmatrix}
\end{equation}
and because $*^2=-1$,
\begin{eqnarray}
B^2-C\cdot A+I&=&0\\
A\cdot B&=&{}^tB\cdot A\\
C\cdot {}^tB&=&B\cdot C.
\end{eqnarray}

On $\diamondsuit$, while the Hodge star $*$ can not be defined, we can
obviously consider the following positive scalar product on the classes of
closed forms modulo $d\varepsilon$, to which the set
$(\alpha^\diamondsuit_k)$ belong:
\begin{eqnarray*}
\notag
  (\alpha^\diamondsuit,\beta^\diamondsuit)&:=&
\left(A(\alpha^\diamondsuit),A(\beta^\diamondsuit)\right)
\\
 &=&
\sum_{\substack{(x,y,x',y')\in\diamondsuit_2\\
 \rho=\rho(x,x'),\,\rho^*=\rho(y,y')}}
\begin{pmatrix}
{\hbox to 0pt{\kern -1.3em t}} \int_{(x\phantom{'},y\phantom{'})}\alpha\\  
\int_{(y\phantom{'},x')}\alpha\\ 
\int_{(x',y')}\alpha\\ \int_{(y',x\phantom{'})}\alpha
\end{pmatrix}
\cdot
\begin{pmatrix}
\scriptscriptstyle +\rho+\rho^*&\scriptscriptstyle+\rho-\rho^*&
\scriptscriptstyle -\rho-\rho^*&\scriptscriptstyle-\rho+\rho^*\\
\scriptscriptstyle +\rho-\rho^*&\scriptscriptstyle+\rho+\rho^*&
\scriptscriptstyle -\rho+\rho^*&\scriptscriptstyle-\rho-\rho^*\\
\scriptscriptstyle -\rho-\rho^*&\scriptscriptstyle-\rho+\rho^*&
\scriptscriptstyle +\rho+\rho^*&\scriptscriptstyle+\rho-\rho^*\\
\scriptscriptstyle -\rho+\rho^*&\scriptscriptstyle-\rho-\rho^*&
\scriptscriptstyle +\rho-\rho^*&\scriptscriptstyle+\rho+\rho^*
\end{pmatrix}
\cdot
\begin{pmatrix}
\int_{(x\phantom{'},y\phantom{'})}\bar\beta\\  
\int_{(y\phantom{'},x')}\bar\beta\\ 
\int_{(x',y')}\bar\beta\\ \int_{(y',x\phantom{'})}\bar\beta
\end{pmatrix}.
\notag  \label{eq:scalarProdDiam} 
\end{eqnarray*}
and it yields
 \begin{equation}
     (\alpha^\diamondsuit_{k},\alpha^\diamondsuit_{\ell})_{k,\ell}
   =
     \begin{pmatrix}
       A_\Gamma+A_{\Gamma^*} & 
{}^tB_{\Gamma\Gamma^*}+{}^tB_{\Gamma^*\Gamma}  \\
      B_{\Gamma\Gamma^*}+B_{\Gamma^*\Gamma} & C_\Gamma+C_{\Gamma^*}
     \end{pmatrix},
     \label{eq:GammaMatrixDiam}
 \end{equation}  
 which, in general, can not be understood as the periods of a set of forms on
 $\diamondsuit$ along the basis $\aleph$.

Let's decompose the space of harmonic forms into two orthogonal supplements,
   \begin{equation}
     \label{eq:parallelPerp}
     \mathcal{H}^1(\Lambda)=
\mathcal{H}^1_\parallel\oplus^\perp \mathcal{H}^1_\perp
   \end{equation}
where the first vector space are the harmonic forms whose holonomies on one
graph are equal to their holonomies on the dual, that is to say
\begin{equation}
  \label{eq:parallel}
  \mathcal{H}^1_\parallel:=\text{Vect }(\alpha_k+\alpha_{k+g}, \;
1\leq k\leq g \text{ ~or~ } 2g< k\leq 3g).
\end{equation}
Definition \eqref{eq:dualBasisDef} and Eq.~\eqref{eq:dualityEta} imply that
\begin{equation}
  \label{eq:perp}
  \mathcal{H}^1_\perp=\text{Vect }(*\alpha_k-*\alpha_{k+g}, \;
1\leq k\leq g \text{ ~or~ } 2g< k\leq 3g).
\end{equation}
These elements in the basis $(\alpha_k+\alpha_{k+g}, \; ;\;
\alpha_k-\alpha_{k+g})$ for $1\leq k\leq g$ and $2g< k\leq 3g$,
are represented by the following invertible matrix:
\begin{equation}
  \label{eq:matrix+-}
  \begin{pmatrix}
    I&0&{}^tB_{\Gamma\Gamma^*}-{}^tB_{\Gamma^*\Gamma}&A_\Gamma-A_{\Gamma^*}\\
    0&I&C_\Gamma-C_{\Gamma^*}&B_{\Gamma\Gamma^*}-B_{\Gamma^*\Gamma}\\
0&0&{}^tB_{\Gamma\Gamma^*}+{}^tB_{\Gamma^*\Gamma}&A_\Gamma+A_{\Gamma^*}\\
0&0&C_\Gamma+C_{\Gamma^*}&B_{\Gamma\Gamma^*}+B_{\Gamma^*\Gamma}
  \end{pmatrix}.
\end{equation}
It implies in particular that the lower right $g\times g$ block is
invertible, therefore so is Eq.~\eqref{eq:GammaMatrixDiam}.

\subsection{Period matrix} \label{sec:periodMatrix}
 \begin{propo} \label{prop:Period}
   The matrix $\Pi=C^{-1}\cdot(i-B)$ is the \emph{period
     matrix}\index{period matrix} of the basis of holomorphic forms
 \begin{equation}
     \zeta_{k}:=(i-*)\sum_{\ell=1}^{2g}C^{-1}_{k,\ell}\;\alpha_{\ell+2g}
     \label{eq:zetaDef}
 \end{equation}  
 in the canonical dissection $\aleph$, that is to say
 \begin{equation}
     \oint_{\aleph_{k}}\zeta_{\ell}=
     \begin{cases}
         \delta_{k,\ell} & \text{~ for ~} 1\leq 
     k\leq 2g, \\
     \Pi_{k-2g,\ell} & \text{~ for ~} 2g< 
         k\leq 4g,
     \end{cases}
     \label{eq:zetaDual}
 \end{equation}
 and $\Pi$ is symmetric, with a positive definite imaginary part.
 \end{propo}
 The proof is essentially the same as in the continuous case~\cite{FK} and we
 include it for completeness.

 \begin{demo}{\ref{prop:Period}}
   Let $\omega_j:=\alpha_j+i*\alpha_j$ for $1\leq j\leq 4g$. These
   holomorphic forms fulfill
   \begin{eqnarray}
     \label{eq:omegaj}
     P_{k,j}:=\tfrac{1}{2}(\omega_k,\omega_j)&=&
(\alpha_k,\alpha_j)+i\,(\alpha_k,-*\,\alpha_j)\\
&=&
     \begin{cases}
        -i\,\int_{\aleph_{j+2g}} \omega_k,& 1\leq j\leq 2g,\\
         i\,\int_{\aleph_{j-2g}} \omega_k,& 2g< j\leq 4g.
     \end{cases}
   \end{eqnarray}
   $P$ is the period matrix of the forms $(\omega)$ in the homology
   basis $\aleph$.  The first $2g$ forms $(\omega_j)_{1\leq j\leq 2g}$
   are a basis of holomorphic forms. It has the right dimension and
   they are linearly independent:
\begin{eqnarray}
  \sum_{j=1}^{2g}(\lambda_j+i\mu_j)(\alpha_j+i*\,\alpha_j)&=&
 \sum_{j=1}^{2g}\left((\lambda_j+\sum_{k=1}^{2g}\mu_k\,
 B_{j,k})\,\alpha_j
+\sum_{k=1}^{2g}\mu_k\, C_{j,k}\,\alpha_{2g+j}\right)\notag\\
&&\; +i\,\sum_{j=1}^{2g}\left((\mu_j+\sum_{k=1}^{2g}\lambda_k\,
 B_{j,k})\,\alpha_j+\sum_{k=1}^{2g}\lambda_k\, C_{j,k}\,\alpha_{2g+j}\right)
  \label{eq:linomega}
\end{eqnarray}
is null, for $\lambda,\mu\in\mathbf{R}$ only when $\lambda=\mu=0$ because $C$
is positive definite. Similarly for the last $2g$ forms. The change of basis
$i\,C^{-1}$ on them provides the basis of holomorphic forms $(\zeta)$. The
last $2g$ rows of $P$ is the $2g\times 4g$ matrix $(B-i\,I,C)$ hence the
periods of $(\zeta)$ in $\aleph$ are given by $(I,\Pi)$.
\end{demo}

The first identity in~Eq.\eqref{eq:zetaDual} uniquely defines the basis
$\zeta$ and a holomorphic $1$-form is completely determined by whether its
periods on the first $2g$ cycles of $\aleph$, or their real parts on the
whole set.

Notice that because $C$ is $g\times g$ block diagonal and $B$ is
anti-diagonal, $\Pi$ is decomposed into four $g\times g$ blocks, the two
diagonal matrices form $i\, C^{-1}$ and are pure imaginary, the other two
form $-C^{-1}\cdot B$ and are real.
 \begin{equation}
   \label{eq:PiBlocks}
   \Pi=
   \begin{pmatrix}
     \Pi_{i*}&\Pi_{r}\\
\Pi_{r*}&\Pi_{i}
   \end{pmatrix}
=
   \begin{pmatrix}
     i\,C_{\Gamma^*}^{-1}&-C_{\Gamma^*}^{-1}\cdot B_{\Gamma^*,\Gamma}\\
-C_{\Gamma}^{-1}\cdot B_{\Gamma,\Gamma^*}&i\,C_{\Gamma}^{-1}
   \end{pmatrix}.
 \end{equation}
 
 Therefore the holomorphic forms $\zeta_k$ are real on one graph and pure
 imaginary on its dual,
 \begin{eqnarray}
   \label{eq:zetaRI}
   1\leq k\leq g&\Rightarrow& \zeta_k\in
   C_{\mathbf{R}}^1(\Gamma)\oplus i\,C_{\mathbf{R}}^1(\Gamma^*)
\\
   g< k\leq 2g&\Rightarrow& \zeta_k\in
   C_{\mathbf{R}}^1(\Gamma^*)\oplus i\,C_{\mathbf{R}}^1(\Gamma).
\notag
 \end{eqnarray}
 
 We will call
 \begin{equation}
   \label{eq:PiGamma}
   \Pi_\Gamma=\Pi_r+\Pi_{i*}
 \end{equation}
 \emph{the period matrix on the graph $\Gamma$} the sum of the real periods
 of $\zeta_{k}$, $1\leq k\leq g$, on $\Gamma$, with the associated pure
 imaginary periods on the dual $\Gamma^{*}$, and similarly for $\zeta_{k}$,
 $g< k\leq 2g$, the period matrix on $\Gamma^{*}$.
 
 It is natural to ask how close $\Pi_{\Gamma}$ and $\Pi_{\Gamma^*}$ are from
 one another, and whether their mean can be given an
 interpretation. Criticality~\cite{M,M01}  answers partially the issue:

 \begin{theo}\label{th:critic}
   In the genus one critical case, the period matrices $\Pi_{\Gamma}$ and
   $\Pi_{\Gamma^*}$ are equal to the period matrix $\Pi_\Sigma$ of the
   underlying surface $\Sigma$. For higher genus, given a refining sequence
   $(\diamondsuit^k,\rho_k)$ of critical maps of $\Sigma$, the discrete
   period matrices $\Pi_{\Gamma^k}$ and $\Pi_{\Gamma^{*k}}$ converge to the
   period matrix $\Pi_\Sigma$.
 \end{theo}
 \begin{demo}{\ref{th:critic}}
   The genus one case is postponed to Sec.~\ref{sec:GenusOne}. The continuous
   limit comes from techniques in~\cite{M,M01}, developed
   in~\cite{M0206041} which prove that, given a refining sequence of critical
   maps, any holomorphic function can be approximated by a sequence of
   discrete holomorphic functions. Taking the real parts, this implies as
   well that any harmonic function can be approximated by discrete harmonic
   functions. In particular, the discrete solutions $f_k$ to a Dirichlet or
   Neumann problem on a simply connected set converge to the continuous
   solution $f$ because the latter can be approximated by discrete harmonic
   functions $g_k$ and the difference $f_k-g_k$ being harmonic and small on
   the boundary, converge to zero. In particular, each form in the basis
   $(\alpha^\diamondsuit_\ell)$, provides a solution to the Neumann problem
   Eq.~\eqref{eq:dualBasisDefDiam} and a similar procedure, detailed
   afterwards, define a converging sequence of forms
   $\zeta^\diamondsuit_\ell$, yielding the result.
 \end{demo}
 
 We can try to replicate the work done on $\Lambda$ on the graph
 $\diamondsuit$.  A problem is that $A_{\Gamma}+A_{\Gamma}$ and
 $C_{\Gamma}+C_{\Gamma}$ need not be positive definite. Moreover, the Hodge
 star $*$ doesn't preserve the space $(A(\alpha^{\Diam}_k))$ of harmonic
 forms with equal holonomies on the graph and on its dual, so we can not
 define the analogue of $\alpha+i\,*\alpha$ on $\diamondsuit$. We first
 investigate what happens when we can partially define these analogues:
 
 Assume that for $2g<k\leq 3g$, the holonomies of $*\alpha_k$ on $\Gamma$ are
 equal to the holonomies of $*\alpha_{k+g}$ on $\Gamma^*$, that is to say
 $C_\Gamma=C_{\Gamma^*}=:\tfrac{1}{2}C_\diamondsuit$ and
 $D_{\Gamma\Gamma^*}=D_{\Gamma^*\Gamma}=:\frac{1}{2}D_\diamondsuit$. It
 implies that the transposes fulfill
 $B_{\Gamma\Gamma^*}=B_{\Gamma^*\Gamma}=:\tfrac{1}{2}B_\diamondsuit$ as well.
 We can then define $\beta^{\Diam}_{k-g}\in Z^1(\diamondsuit)$ such that
 $A(\beta^{\Diam}_{k-g})=*\alpha_{k+g}$, uniquely up to $d\varepsilon$. The
 last $g$ columns ${}^t(B_\diamondsuit,C_\diamondsuit)$ of the matrix of
 scalar product Eq.~\eqref{eq:GammaMatrixDiam} are related to their periods
 in the homology basis $\aleph^{\Diam}$ in a way similar to
 Eq.~\eqref{eq:GammaMatrix}. By the same reasoning as before, the forms
\begin{equation}
  \label{eq:zetaDiam}
  \zeta^\diamondsuit_k=\sum_{\ell=1}^g{C_\diamondsuit^{-1}}_{k,\ell}
\left(\alpha^{\Diam}_{\ell+g}-i\beta^{\Diam}_{\ell+g}\right), \;\;
1\leq k\leq g
\end{equation}
verify $A(\zeta^\diamondsuit_k)=\frac{\zeta_k+\zeta_{k+g}}{2}$ and have
periods on $\aleph^{\Diam}$ given by the identity for the first $g$ cycles
and the following $g\times g$ matrix, mean of the period matrices on the
graph and on its dual:
\begin{equation}
\Pi^{\Diam}=C_\diamondsuit^{-1}(i-B_\diamondsuit)=
\frac{\Pi_{\Gamma}+\Pi_{\Gamma^*}}{2}.
     \label{eq:PeriodMatrixDiamLambda}
 \end{equation}
 
 The same reasoning applies when the periods of the forms $*\alpha_k$ on the
 graph and on its dual are not equal but close to one another. In the context
 of refining sequences, we said that the basis $(\alpha^\diamondsuit_\ell)$,
 converges to the continuous basis of harmonic forms defined by the same
 Neumann problem Eq.~\eqref{eq:dualBasisDefDiam}. Therefore
 \begin{equation}
   \label{eq:CGammaCGammaS}
   C_\Gamma-C_{\Gamma^*}=o(1), \qquad
 B_{\Gamma\Gamma^*}-B_{\Gamma^*\Gamma}=o(1).
 \end{equation}
 A harmonic form $ \nu_{k+g}=o(1) $ on $\Gamma^*$ can be added to $
 *\alpha_{k+g} $ such that there exists $ \beta^{\Diam}_{k-g}\in
 Z^1(\diamondsuit) $ with $ A(\beta^{\Diam}_{k-g})=*\alpha_{k+g}+\nu_{k+g}$,
 yielding forms $ \zeta^\diamondsuit_k$, verifying $
 A(\zeta^\diamondsuit_k)=\tfrac{1}{2}(\zeta_k+\zeta_{k+g})+o(1) $ and whose
 period matrix is $ \Pi^{\Diam}+o(1)$. Since the periods of $\alpha_k$
 converge to the same periods as its continuous limit, this period matrix
 converges to the period matrix $\Pi_\Sigma$ of the surface. Which is the
 claim of Th.~\ref{th:critic}.

 In the paper~\cite{CSMcC}, R.~Costa-Santos and B.~McCoy define a period
 matrix on a special cellular decomposition $\Gamma$ of a surface by squares.
 They don't consider the dual graph $\Gamma^{*}$.  Their period matrix is
 equal to one of the two diagonal blocks of the double period matrix we
 construct in this case. They don't have to consider the off-diagonal blocks
 because the problem is so symmetric that their period matrix is pure
 imaginary.

\subsection{Genus one case} \label{sec:GenusOne}
Criticality solves partially the problem of having two different $g\times g$
period matrices instead of one since they converge to one another in a
refining sequence. However, on a genus one critical torus, the situation is
simpler: The overall curvature is null and a critical map is everywhere flat.
Therefore the cellular decomposition is the quotient of a periodic cellular
decomposition of the plane by two independant periods. They can be normalized
to $(1,\tau)$.  The continuous period matrix is the $1\times 1$-matrix
$\tau$. A basis of the two dimensional holomorphic $1$-forms is given by the
real and imaginary parts of $dZ$ on $\Gamma$ and $\Gamma^{*}$ respectively,
and the reverse. The discrete period matrix is the $2\times 2$ matrix $
\begin{pmatrix}
    \text{Im }\tau& \text{Re }\tau\\
    \text{Re }\tau& \text{Im }\tau
\end{pmatrix}
$ and the period matrices on the graph and on its dual are both equal to the
continuous one.

For illustration purposes, the whole construction, of a basis of
harmonic forms, then projected onto a basis of holomorphic forms,
yielding the period matrix, can be checked explicitely on the critical
maps of the genus $1$ torus decomposed by square or
triangular/hexagonal lattices:

Consider the critical square (rectangular) lattice decomposition of a
torus $\Diam=(\mathbb{Z}e^{i\,\theta}+\mathbb{Z}e^{-i\,\theta})/
(2p\,e^{i\,\theta}+2q\,e^{-i\,\theta})$, with horizontal parameter
$\rho=\tan\theta$ and vertical parameter its inverse. Its modulus is 
$\tau=\frac{q}{p}e^{2\,i\,\theta}$. The two dual
graphs $\Gamma$ and $\Gamma^{*}$ are isomorphic.  An explicit harmonic
form $\alpha^{\Gamma}_{1}$ is given by the constant $1/2p$ on
horizontal and downwards edges of the graph $\Gamma$ and $0$ on all
the other edges.  Its holonomies are $1$ and $0$ on the $p$, resp. 
$q$ cycles.  Considering $1/2q$ and the dual graph, we construct in
the same fashion $\alpha^{\Gamma}_{2}, \alpha^{\Gamma^{*}}_{1},
\alpha^{\Gamma^{*}}_{2}$.  The matrix of inner products is
\begin{equation}
    (\alpha_{k},\alpha_{\ell})_{k,\ell}
    =  \frac{1}{\sin 2\theta}
    \begin{pmatrix}
        \frac{q}{p} &&& \cos 2\theta  \\
        &\frac{q}{p} & \cos 2\theta   \\
    & \cos 2\theta&     \frac{p}{q}\\
    \cos 2\theta&&&     \frac{p}{q}
   \end{pmatrix}
    \label{eq:GammaMatrixSq}
\end{equation} 
using $\frac{\rho+1/\rho}{2}=1/\sin 2\theta$ and 
$\frac{\rho-1/\rho}{2}=-1/\tan 2\theta$ so that the period matrix is
\begin{equation}
    \Pi=\frac{q}{p}
    \begin{pmatrix}
        i\,\sin 2\theta&\cos 2\theta\\
        \cos 2\theta&i\,\sin 2\theta
    \end{pmatrix}.
    \label{eq:PeriodMatrixSq}
\end{equation}
Therefore there exists a holomorphic form which has the same periods
on the graph and on its dual, it is the average of the two half forms
of Eq.~\eqref{eq:zetaDual} and its periods are
$(1,\frac{q}{p}e^{2\,i\,\theta})$ along the $p$, resp.  $q$ cycles,
yielding the continuous modulus.  This holomorphic form is simply
the normalized fundamental form $\frac{d Z}{p e^{-i\,\theta}}$.
 
In the critical triangular/hexagonal lattice, we just point out to the
necessary check by concentrating on a tile of the torus, composed of
two triangles, pointing up and down respectively. We show that there 
exists an explicit holomorphic form which has the same shift on the 
graph and on its dual, along this tile.  Let $\rho_{-},
\rho_{\backslash}$ and $\rho_{/}$ the three parameters around a given
triangle.  Criticality occurs when $\rho_{-}\,\rho_{\backslash}+
\rho_{\backslash}\,\rho_{/}+ \rho_{/}\,\rho_{-}=1$.  The form which is
$1$ on the rightwards and South-West edges and $0$ elsewhere is
harmonic on the triangular lattice.  Its pure imaginary companion on
the dual hexagonal lattice exhibits a shift by $i\,\rho_{\backslash}$
in the horizontal direction and $i\,(\rho_{\backslash}+\rho_{-})$ in
the North-East direction along the tile.  Dually, on the hexagonal
lattice, the form which is $\rho_{\backslash}\,\rho_{-}$ along the
North-East and downwards edges and $1-\rho_{\backslash}\,\rho_{-}$
along the South-East edges, is a harmonic form. Its shift in the
horizontal direction is $1$, in the North-East direction $0$, and its
pure imaginary companion on the triangular lattice exhibits a shift by
$i\,\rho_{\backslash}$ in the horizontal direction and
$i\,(\rho_{\backslash}+\rho_{-})$ in the North-East direction along
the tile as before.  Hence their sum is a holomorphic form with equal
holonomies on the triangular and hexagonal graphs and the period
matrix it computes is the same as the continuous one.  This simply
amounts to pointing out that the fundamental form $dz$ can be explicitely
expressed in terms of the discrete conformal data.

\section{Criticality and integrable system}
This theory can be viewed as the simplest (it is linear) of a series
of integrable theories~\cite{M0402097}. We will present its quadratic counterpart,
which leads to another version of discrete analytic functions, based
on circle patterns. Along the way, we will see how discrete
exponentials and discrete polynomials emerge due to integrable systems
theory pieces of technology named the Bäcklund or Darboux transform~\cite{M0402097}.

\subsection{Criticality}
Until now, everything has been purely combinatorial, there was no
reference to an underlying geometry and no continuous limit.
Criticality is what links combinatorics and geometry, and what gives a
meaning to approximation theorems.

\begin{defi}\index{criticality}
  A discrete conformal map $(\diamondsuit,\rho)$ is critical if there
  exists a discrete holomorphic map $Z$ such that the quadrilateral
  faces $\diamondsuit_2$ can be simultaneously embedded into \textbf{rhombi} in
  the complex plane. 
\end{defi}

Because of the Gauss-Bonnet theorem, it is not possible to globally
embed a compact surface into the plane, therefore we allow for an
atlas of local critical maps with a finite number of fixed local conic
singularities. When a continuous limit is taken, their number, angle
and position should not change, and the theorem of isolated
singularities helps us wipe them out as inessential.

It is a simple calculation to check that if $Z$ is a critical map, any
discrete holomorphic function $f\in\Omega(\diamondsuit)$ gives rise,
through \eqref{eq:wedge01} to a holomorphic $1$-form $f dZ$.

For a holomorphic function $f$, the equality $f\,dZ\equiv 0$ is
equivalent to $f=\lambda\,\varepsilon$ for some $\lambda\in
\mathbb{C}$ with $\varepsilon$ the {biconstant}
$\varepsilon(\Gamma)=+1$, $\varepsilon(\Gamma^*)=-1$.

Following Duffin~\cite{Duf, Duf68}, we introduce the
\begin{defi}\label{def:fp}
  For a holomorphic function $f$, define on a flat simply connected map $U$
  the holomorphic functions $f^{\dag}$, the \emph{dual} of $f$, and $f'$, the
  \emph{derivative} of $f$, by the following formulae:\index{discrete derivative}
  \begin{equation}
f^{\dag}(z):=\varepsilon(z)\,\bar f(z),\label{eq:fdag}
\end{equation}
where $\bar f$ denotes the
  complex conjugate, $\varepsilon=\pm 1$ is the biconstant, and
  \begin{equation}
f'(z):=\frac4{\delta^2}\left( \int_{O}^z f^{\dag}
    dZ\right)^{\dag}+\lambda\,\varepsilon,\label{eq:deffp}
  \end{equation}
defined up to $\varepsilon$.
\end{defi}

It is an immediate  calculation~\cite{M01} to check the following 
\begin{propo}\label{prop:fp}
  The derivative $f'$ fulfills
\begin{equation}
  \label{eq:dffpdz}
  d\, f = f'\, dZ.
\end{equation}
\end{propo}

\subsection{$\bar\partial$ operator}\index{$\bar\partial$ operator}
For holomorphic or anti-holomorphic functions, $df$ is, locally on
each pair of dual diagonals, proportional to $dZ$, resp. $d\bar Z$, we
define $\partial$ and $\bar\partial$ operators (not to be confused
with the boundary operator on chains) that decompose the exterior
derivative into holomorphic and anti-holomorphic parts yielding
$$df\wedge\overline{df} = \left(|\partial f|^2 +|\bar\partial
  f|^2\right) dZ\wedge d\bar Z$$
where the derivatives naturally live
on faces:

In the continuous theory, for any derivable function $f$ on the
complex plane, the derivatives $\partial =
\frac{d\phantom{x}}{dx}+i\,\frac{d\phantom{y}}{dy}$ and $\bar\partial
= \frac{d\phantom{x}}{dx}-i\,\frac{d\phantom{y}}{dy}$ with respect to
$z=x+i\,y$ and $\bar z=x-i\,y$ yield
$$f(z+z_0)=f(z_0)+z(\partial f)(z_0) +\bar z(\bar\partial f)(z_0)
+o(|z|).$$

These derivatives can be seen as a limit of a contour integral over
a small loop $\gamma$ around $z_0$:
$$(\partial f)(z_0)=\lim_{\gamma\to z_0}\frac
i{2\mathcal{A}(\gamma)}{\displaystyle\oint_{\gamma}f d\bar z},
\qquad\qquad
(\bar\partial f)(z_0)=-\lim_{\gamma\to z_0}\frac
i{2\mathcal{A}(\gamma)}{\displaystyle\oint_{\gamma}f d Z},
$$
which leads to the following definitions in the
discrete setup:
$$\begin{array}{rcl}
\partial:C^0(\Diam)&\to&C^2(\Diam)\\
f&\mapsto&\partial f = \bigl[(x,y,x',y')\mapsto \frac
i{2\mathcal{A}(x,y,x',y')}{\displaystyle\oint\limits_{(x,y,x',y')}f d\bar
  Z}\bigr],\\ 
\\
\bar\partial:C^0(\Diam)&\to&C^2(\Diam)\\
f&\mapsto&\bar\partial f = \bigl[(x,y,x',y')\mapsto -\frac
i{2\mathcal{A}(x,y,x',y')}{\displaystyle\oint\limits_{(x,y,x',y')}f
  dZ}\bigr].\\ 
\end{array}
$$
A holomorphic function $f$ verifies $\bar\partial f\equiv 0$ and (with
$Z(u)$ noted simply $u$)
$$\partial f(x,y,x',y') = \frac{f(y')-f(y)}{y'-y} =
\frac{f(x')-f(x)}{x'-x}.$$
 Notice that the statement
$f=\left(\int\partial f\,dz\right)$ has no meaning, $\partial$ is not
a derivation endomorphism in the space of functions on the vertices of
the double.

On the other hand, these differential operators can be extended
(see~\cite{Ken02}) into operators (the Kasteleyn operator)
$\partial_{20},\bar\partial_{20}:C^2(\Diam)\to C^0(\Diam)$ simply by
transposition, $\partial_{20} = -^t\partial_{02}$, leading to
endomorphisms of $C^0(\Diam)\oplus C^2(\Diam)$. They are such that
their composition, restricted to the vertices $\Diam_0$, gives back
the laplacian:
$$\Delta =\tfrac12\left(
 \partial\circ\bar\partial + \bar\partial\circ\partial\right).$$

Furthermore, the double derivative $\partial_{20}\circ\partial_{02}$
is a well defined endomorphism of $C^0(\Diam)$.

\goodbreak
\subsection{Discrete exponential}\label{sec:ExpDef}
\begin{defi}\index{discrete exponential}
  For a constant $\lambda\in\mathbb{C}$, the discrete exponential
  $\exp({:}\lambda{:}\,Z)$ is the solution of
\begin{eqnarray}
\exp({:}\lambda{:}\,O)&=&1\label{eq:Exp}\notag\\
d\,\exp({:}\lambda{:}\,Z)&=&\lambda\,\exp({:}\lambda{:}\,Z)\,dZ.\label{eq:dExp}
\end{eqnarray}
We define its derivatives with respect to the continuous parameter $\lambda$:
\begin{equation}
  \label{eq:ZkExp}
  Z^{{:}k{:}}\exp({:}\lambda{:}\,Z)\, := \, \frac{\partial^k~}{\partial\lambda^k}\,\exp({:}\lambda{:}\,Z).
\end{equation}
\end{defi}

The discrete exponential on the square lattice was defined by
Lelong-Ferrand~\cite{LF}, generalized in~\cite{M0111043} and studied
independently in~\cite{Ken02,Ken02I}.  For
$\abs{\lambda}\not=2/\delta$, an immediate check shows that it is a
rational fraction in $\lambda$ at every point: For the vertex
$x=\sum\delta\,e^{i\,\theta_k}$,
\begin{equation}
  \label{eq:ExpRatFrac}
  \exp({:}\lambda{:}\,x)=\prod_{k}
\frac{1+\frac{\lambda\delta}{2}e^{i\,\theta_k}}{
1-\frac{\lambda\delta}{2}e^{i\,\theta_k}}
\end{equation}
where $(\theta_k)$ are the angles defining
$(\delta\,e^{i\,\theta_k})$, the set of ($Z$-images of)
$\diamondsuit$-edges between $x$ and the origin.  Because the map is
critical, Eq.~\eqref{eq:ExpRatFrac} only depends on the end points
$(O,x)$. It is a generalization of a well known formula, in a slightly
better version,
\begin{equation}
  \label{eq:ExpLimn}
  \exp(\lambda\, x) = \( 1+\frac{\lambda\, x}{n}\)^n +O(\frac{\lambda^2\,
    x^2}{n}) =  
\( \frac{1+\frac{\lambda\, x}{2n}}{1-\frac{\lambda\, x}{2n}}\)^n +O(\frac{\lambda^3\, x^3}{n^2})
\end{equation}
to the case when the path from the origin to the point $x=\sum_1^n
\frac{x}{n}=\sum\delta\,e^{i\,\theta_k}$ is not restricted to straight equal segments
but to a general path of $O(|x|/\delta)$ segments of any directions.

The integration with respect to $\lambda$ gives an interesting analogue of
$Z^{{:}-k{:}}\,\exp({:}\lambda{:}\,Z)$. It is defined up to a globally defined
discrete holomorphic map. One way to fix it is to integrate from a given
$\lambda_0$ of modulus $2/\delta$, which is not a pole of the rational
fraction, along a path that doesn't cross the circle of radius $2/\delta$
again.

\begin{propo}\label{prop:ExpDag}
For point-wise multiplication, at every point
$x\in\diamondsuit_0$, 
\begin{equation}
\exp({:}\lambda{:}\,x)\cdot\exp({:}-\lambda{:}\,x)=1.\label{eq:ExpExp1}
\end{equation}
The specialization at $\lambda=0$ defines monomials:
\begin{equation}
 Z^{{:}k{:}}\, = \, Z^{{:}k{:}}\exp({:}\lambda{:}\,Z)|_{\lambda=0}
\label{eq:ZkExpZk}
\end{equation}
which fulfill $Z^{{:}k{:}}\, = \, k\,\int Z^{{:}k-1{:}}\, dZ$.
The anti-linear duality ${\dag}$ maps exponentials to exponentials:
\begin{equation}
  \label{eq:ExpDag}
  \exp({:}\lambda{:})^{\dag}=\exp({:}\frac{4}{\delta^2\bar\lambda}{:}).
\end{equation}
\end{propo}
In particular, $ \exp({:}\infty{:})=1^{\dag}=\varepsilon$ is the biconstant.
\begin{proof}[Proof \ref{prop:ExpDag}]
The first assertion is immediate.

The derivation of \eqref{eq:dExp} with respect to $\lambda$ yields
\begin{equation}
  \label{eq:dZkExp}
  d\, \frac{\partial^k~}{\partial\lambda^k}\,\exp({:}\lambda{:}\,Z) = \(
\lambda\, \frac{\partial^k~}{\partial\lambda^k}\,\exp({:}\lambda{:}\,Z) 
+ k \, \frac{\partial^{k-1}~}{\partial\lambda^{k-1}}\,\exp({:}\lambda{:}\,Z)
\) dZ
\end{equation}\index{discrete polynomials}
which implies \eqref{eq:ZkExpZk}.

Derivation of $\exp({:}\lambda{:})^{\dag}$ gives, 
\begin{eqnarray}
  \label{eq:ExpDagDeriv}
\(\exp({:}\lambda{:})^{\dag}\)'&=&
\frac{4}{\delta^2}\(\int_O^z \exp({:}\lambda{:})\, dZ\)^{\dag}+\mu\,\varepsilon\\
&=&\frac{4}{\delta^2}\,\(\frac{\exp({:}\lambda{:})-1}{\lambda}\)^{\dag}
+\mu\,\varepsilon
\notag\\
&=&\frac{4}{\delta^2\bar\lambda}\,\exp({:}\lambda{:})^{\dag}+\nu\,\varepsilon
\notag
\end{eqnarray}
with $\mu,\,\nu$ some constants, so that the initial condition
$\exp({:}\lambda{:}O)^{\dag}=1$ at the origin and the difference equation
$d\,\exp({:}\lambda{:})^{\dag}=
\frac{4}{\delta^2\bar\lambda}\,\exp({:}\lambda{:})^{\dag}\,dZ$ yields the
result.
\end{proof}

Note that it is natural to define
$\exp({:}\lambda{:}\,(x-x_0)):=\frac{\exp({:}\lambda{:}\,x)}{\exp({:}\lambda{:}\,x_0)}$
as a function of $x$ with $x_0$ a fixed vertex. It is simply a change of
origin. But apart on a lattice where addition of vertices or multiplication
by an integer can be given a meaning as maps of the lattice, there
is no easy way to generalize this construction to other discrete holomorphic
functions such as $\exp({:}\lambda{:}\,(x+n\,y))$ with $x,y\in\diamondsuit_0$
and $n\in\mathbb{Z}$.

\subsection{Series}\label{sec:Series}\index{discrete series}
The series $\sum_{k=0}^{\infty}\frac{\lambda^{k}\,Z^{:k:}}{k!}$, wherever it
is absolutely convergent, coincide with the rational
fraction~\eqref{eq:ExpRatFrac}: Its value at the origin is $1$ and it
fulfills the defining difference equation~\eqref{eq:dExp}. Using
Eq.~\eqref{eq:ZkExpZk}, a Taylor expansion of $\exp({:}\lambda{:}\,x)$ at
$\lambda=0$ gives back the same result. We are now interested in the rate of
growth of the monomials.

Direct analysis gives an estimate of $Z^{:k:}$:
\begin{propo}\label{prop:ZkCrois}
  For $x\in\diamondsuit$, at a combinatorial distance $d(x,O)$ of the origin,
  and any $k\in\mathbb{N}$,
  \begin{equation}
    \label{eq:ZkCrois}
    \abs{\frac{Z^{:k:}(x)}{k!}}\leq\(\frac{\alpha+1}{\alpha-1}\)^{{d({x},O)}}
\(\alpha\,\frac{\delta}{2}\)^k,
  \end{equation}
for any $\alpha>1$ arbitrarily close to $1$.
\end{propo}

\begin{corol}\label{prop:ExpSeriesConv}
  The series  $\sum_{k=0}^{\infty}\frac{\lambda^{k}\,Z^{:k:}}{k!}$ is
  absolutely convergent for $\abs{\lambda}<\frac{2}{\delta}$.
\end{corol}

\begin{proof}[Proof \ref{prop:ZkCrois}]
  It is proved by double induction, on the degree $k$ and on the
  combinatorial distance to the origin.  
  
  For $k=0$, it is valid for any $x$ since
  $\frac{\alpha+1}{\alpha-1}=1+\frac{2}{\alpha-1}>1$, with equality only at
  the origin.
  
  Consider $x\in\diamondsuit$ a neighbor of the origin,
  $Z(x)=\delta\,e^{i\,\theta}$, then an immediate induction gives for $k\geq
  1$,
  \begin{equation}
    \label{eq:ZkxO}
    \frac{Z^{:k:}(x)}{k!}=2\(\frac{\delta\,e^{i\,\theta}}{2}\)^k
  \end{equation}
  which fulfills the condition Eq.~\eqref{eq:ZkCrois} for any $k\geq1$
  because $\frac{\alpha+1}{\alpha-1}\,\alpha^k>2$. This was done merely for
  illustration purposes since it is sufficient to check that the condition
  holds at the origin, which it obviously does.  
  
  Suppose the condition is satisfied for a vertex $x$ up to degree $k$, and
  for its neighbor $y$, one edge further from the origin, up to degree $k-1$.
  Then,
\begin{equation}
  \label{eq:Zkxy}
  \frac{Z^{:k:}(y)}{k!}=\frac{Z^{:k:}(x)}{k!}+
\frac{Z^{:k-1:}(x)+Z^{:k-1:}(y)}{(k-1)!}\frac{Z(y)-Z(x)}{2}
\end{equation}
 in absolute value fulfills
\begin{eqnarray}
  \abs{\frac{Z^{:k:}(y)}{k!}}&\leq&
  \(\frac{\alpha+1}{\alpha-1}\)^{{d({x},O)}}\(\alpha\,\frac{\delta}{2}\)^{k-1}
\(
\(\alpha\,\frac{\delta}{2}\)+
\({1+\frac{\alpha+1}{\alpha-1}}\)\frac{\delta}{2}\)
\notag
\\
&& =\(\frac{\alpha+1}{\alpha-1}\)^{{d({x},O)}}\(\alpha\,\frac{\delta}{2}\)^{k}
\(1+\frac{2}{\alpha-1}\)  \label{eq:ZkxyAbs}
\\
&& =\(\frac{\alpha+1}{\alpha-1}\)^{{d({y},O)}}\(\alpha\,\frac{\delta}{2}\)^{k},
\notag
\end{eqnarray}
thus proving the condition for $y$ at degree $k$. It follows by induction
that the condition holds at any point and any degree.
\end{proof}

\subsection{Basis} The discrete exponentials form a basis of discrete
  holomorphic functions on a finite critical map: given any set of
  pair-wise different reals $\{\lambda_k\}$ of the right dimension,
  the associated discrete exponentials will form a basis of the space
  of discrete holomorphic functions. See~\cite{M0402097} for the
  formula 
  \begin{equation}
    \label{eq:Fourier}
    f(x)=\frac1{2i\pi}\int_\gamma g(\lambda) \exp(:\lambda: x)\, d\lambda
  \end{equation}
  for a certain fixed contour $\gamma$ in the space of parameters
  $\lambda$, and the definition of $g(\lambda)$ as a fixed contour
  integral in $\diamondsuit$ involving $f$.
  
  The polynomials however don't form a basis in general: the
  combinatorial surface has to fulfill a certain condition called
  ``combinatorial convexity''. A quadrilateral, when traversed from
  one side to its opposite, define a unique chain of quadrilaterals,
  that we call a ``train-track''. The condition we ask is that two
  different train-tracks have different slopes.
  
  On a combinatorially convex set, the discrete polynomials form a
  basis as well.

\subsection{Continuous limit}

In a critical map, where quadrilaterals are mapped to rhombi of side
$\delta$, identifying a vertex $x$ with its image $Z(x)$.

The
combinatorial distance $d_\diamondsuit(x,O)$ is related to the modulus
$|x|$ through
\begin{equation}
  \label{eq:dist}
  d_\diamondsuit(x,O) \frac{\sin\theta_m}4 \leq \frac{|x|}{\delta}\leq d_\diamondsuit(x,O)
\end{equation}
where $\theta_m$ is the minimum of all rhombi angles. When the
rhombi don't flatten, the combinatorial distance and the modulus (over
$\delta$) are equivalent distances.

\begin{lemma}           \label{lem:diam}
  Let $(ABCD)$ be a four sided polygon of the Euclidean plane such
  that its diagonals are orthogonal and the vertices angles are in
  $[\eta, 2\pi-\eta]$ with $\eta>0$. Let $(M,M')$ be a pair of points
  on the polygon. There exists a path on $(ABCD)$ from $M$ to $M'$ of
  minimal length $\ell$. Then
  $$\frac{MM'}\ell\geq\frac{\sin \eta}4.$$
\end{lemma}

\begin{figure}[htbp]
\begin{center}\input{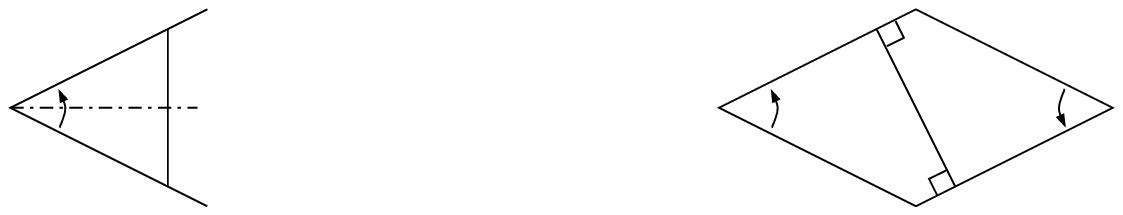}
\end{center}
\caption{The two extremal positions.}   \label{fig:diamlemma}

\end{figure}

It is a straightforward study of a several variables function. If the
two points are on the same side, $MM'=\ell$ and $\sin\eta\leq 1$. If
they are on adjacent sides, the extremal position with $MM'$ fixed is
when the triangle $MM'P$, with $P$ the vertex of $(ABCD)$ between
them, is isocel. The angle in $P$ being less than $\eta$,
$\frac{MM'}\ell\geq{\sin \frac\eta2}>\frac{\sin \eta}2.$ If the points
are on opposite sides, the extremal configuration is given by
Fig.~\ref{fig:diamlemma}.2., where $\frac{MM'}\ell=\frac{\sin
  \eta}4$.$\Box$

A function $f:\diamondsuit_0\to\mathbb{C}$ on the combinatorial
surface can be extended to a function on the image of the critical map
in the complex plane $\hat{f}:U\to\mathbb{C}$ by stating that 
$\hat{f}\(Z(x)\)=f(x)$ for the image of a vertex, and extend it
linearly on the segments $[Z(x),Z(y)]$ image of an edge, and
harmonically inside each rhombus.

\begin{theo}\label{th:limit}\index{continuous limit}
  Let $(\diamondsuit_k)$ a sequence of simply connected critical maps,
  $U$ the non empty intersection of their images in the complex plane
  and a holomorphic function $f:U\to\mathbb{C}$. If the sequence of
  minimum angles are bounded away from $0$ and the sequence of rhombi
  side lengths $(\delta_k)$ converge to $0$, then the function $f$ can
  be approximated by a sequence of discrete holomorphic functions
  $f_n\in\Omega(\diamondsuit_k)$ converging to $f$. The convergence is
  not only pointwise but $C^\infty$ on the intersection of images.
  Conversely a converging sequence of discrete holomorphic functions
  converges to a continuous holomorphic function, in particular the
  discrete polynomials and discrete exponentials with fixed
  parameters.
\end{theo}

\begin{corol}
  On a Riemann surface, any $1$-form can be approximated by a sequence
  of discrete holomorphic $1$-forms on a refining sequence of critical
  maps with fixed conic singularities.
\end{corol}

The proof relies on the convergence of polynomials seen as iterated
primitives of the constant function.
\begin{lemma}\label{th:primitive}
  Given a sequence of discrete holomorphic functions $(f_k)$ on a refining
  sequence of critical maps, converging to a holomorphic function $f$, the
  sequence of primitives $(\int f_k\, dZ)$ converges to $\int f(z)\, dz$.
  Moreover, in the compact case, if the convergence of the functions is of
  order $O(\delta_k^2)$, it stays this way for the primitives.
\end{lemma}

\begin{demo}{\ref{th:primitive}}
  Suppose that we are given a sequence of flat vertices
  $O_k\in\diamondsuit_k$ where the face containing the fixed flat
  origin $O\in U$ is adjacent to $O_k$.  For a given integer $k$, let
  $\widehat F_k$ the continuous piecewise harmonic extension of the
  discrete primitive $\int_{O_k} f_k\, dZ$ to $U$. We want to prove
  that for any $x\in U$, the following sequence tends to zero
\begin{equation}
  \label{eq:uniformF}
\(\left|(\widehat F_k(x)-\widehat F_k(O))-\int_O^x f(z)\, 
dz\right|\)_{k\in\mathbb{N}}.
\end{equation}

For each integer $k$ consider a vertex $x_k\in\diamondsuit_0$ on the
boundary of the face of $\diamondsuit_2$ containing $x$.

We decompose the difference~\eqref{eq:uniformF} into three parts, inside the
face containing the origin $O$ and its neighbor $O_k$, similarly for $x$ and
$x_k$, and purely along the edges of the graph $\diamondsuit^k$ itself.
\begin{eqnarray}
|\bigl(\widehat F_k(x)&-&\widehat F_k(O)\bigr)-\int_O^x f(z)\, dz| \;=\;\notag\\&&
\left|(\widehat F_k(x)-\widehat F_k(x_k))+\int_{O_k}^{x_k}f_k\, dZ+
(\widehat F_k(O)-\widehat F_k(O_k))-
\int_O^x f(z)\, dz\right|\notag\\
&\leq&
\left|\widehat F_k(x)-\widehat F_k(x_k)-\int_{x_k}^x f(z)\, dz\right|+
\left|\int_{O_k}^{x_k}f_k\, dZ-\int_{O_k}^{x_k} f(z)\, dz\right|+\notag\\
&&\left|\widehat F_k(O_k)-\widehat F_k(O)-\int_{O}^{O_k} f(z)\, dz\right|.  
\label{eq:F}
\end{eqnarray}
On the face of $\diamondsuit$ containing $x$, the primitive $\xi\mapsto
\int_{x_k}^\xi f(z)\, dz$ is a holomorphic, hence harmonic function as
well as $\xi\mapsto \widehat F_k(\xi)$.  By the maximum principle, the
harmonic function $\xi\mapsto \widehat F_k(\xi)-\widehat
F_k(x_k)-\int_{x_k}^\xi f(z)\, dz$ reaches its maximum on that face,
along its boundary. The difference of the discrete primitive along the
boundary edge $(x_k,y)\in\diamondsuit_1$ at the point
$\xi=(1-\lambda)x_k+\lambda\, y$ is equal by definition to
\begin{equation}
\widehat F_k((1-\lambda)x_k+\lambda\, y)-\widehat
F_k(x_k)=\lambda(y-x_k)\frac{f_k(x_k)+f_k(y)}2.
\label{eq:Fdifint}
\end{equation}
The holomorphic $f$ is differentiable with a bounded derivative on
$U$, so averaging the first order expansions at $x_k$ and $y$, we get
\begin{eqnarray}
  \int_{x_k}^{\xi} f(z)\, dz&=&
\lambda(y-x_k)\frac{f(x_k)+f(y)}2+(y-x_k)^2\frac{\lambda^2f'(x_k)+(1-\lambda)^2f'(y)}4
\notag\\
&&\qquad\qquad\qquad\qquad+o\left((\xi-x_k)^3\right)+o\left((\xi-y)^3\right)\notag\\
&=&\lambda(y-x_k)\frac{f(x_k)+f(y)}2+O(\delta_k^2)\notag\\
  \label{eq:fint}
\end{eqnarray}
therefore
\begin{equation}
  \label{eq:Fx}
  |\widehat F_k(x)-\widehat F_k(x_k)-\int_{x_k}^x f(z)\, dz|=O(\delta_k^2).
\end{equation}
  Similarly for the term around the origin.
  
  By definition of $\widehat f_k$, the $1$-form $\widehat f_k(z)\, dz$ along
  edges of the graph $\diamondsuit$ is equal to the discrete form $f_k dZ$ so
  that $\int_{O_k}^{x_k}f_k\, dZ=\int_{O_k}^{x_k}\widehat f_k(z)\, dz$ on a
  path along $\diamondsuit$ edges.
  Therefore the difference
\begin{equation}
  \label{eq:Fint}
\left|\int_{O_k}^{x_k}f_k\, dZ-\int_{O_k}^{x_k} f(z)\, dz\right|\leq
\int_{O_k}^{x_k}\left|\(\widehat f_k(z)-f(z)\)\, dz\right|
\end{equation}
is of the same order as the difference $\left|f_k(z)-f(z)\right|$
times the length $\ell(\gamma_k)$ of a path on $\diamondsuit_k$ from
$O_k$ to $x_k$. This length is bounded as $\ell(\gamma_k)\leq \frac
4{\sin\theta_m}|x_k-O_k|$.  Since we are interested in the compact
case, this length is bounded uniformly and the
difference~\eqref{eq:Fint} is of the same order as the point-wise
difference.  We conclude that the sequence of discrete primitives
converges to the continuous primitive and if the limit for the
functions was of order $O(\delta^2)$, it remains of that order.
\end{demo}

The discrete polynomials of degree less than three agree point-wise with their
continuous counterpart, $Z^{:2:}(x)=Z(x)^2$.

A simple induction then gives the following
\begin{corol}\label{cor:poly}
  The discrete polynomials converge to the continuous ones, the limit is
  of order $O(\delta_k^2)$.
\end{corol}

Which implies the main theorem:

\begin{demo}{\ref{th:limit}}
  On the simply connected compact set $U$, a holomorphic function $f$ can be
  written, in a local map $z$ as a series,
  \begin{equation}
    \label{eq:series}
    f(z)=\sum_{k\in\mathbb{N}} a_k z^k.
  \end{equation}
  Therefore, by a diagonal procedure, there exists an increasing
  integer sequence $\(N(n)\)_{n\in\mathbb{N}}$ such that the sequence
  of discrete holomorphic polynomials converge to the continuous
  series and the convergence is $C^\infty$.
\begin{equation}
  \label{eq:seriesPoly}
  \(\sum_{k=0}^{N(n)} a_k Z^{:k:}\)_{n\in\mathbb{N}}\to f.
\end{equation}
\end{demo}

\subsection{Cross-ratio preserving maps}

Once the isometry $Z$ is chosen, holomorphicity of a function $f$ can be written
on a quadrilateral $(x,y,x',y')\in\diamondsuit_2$, writing $x=Z(x)$ for a vertex
$x\in\diamondsuit_0$, as
\begin{equation}
  \label{eq:holoZ}
  \frac{f(y')-f(y)}{f(x')-f(x)}=
  \frac{y'-y}{x'-x}
\end{equation}
and $f$ is understood as a \emph{diagonal ratio preserving map}, and
each value at a corner vertex can be linearly solved in terms of the
three others.

A quadratic version is given by the \emph{cross-ratio preserving
  maps}: A function $f$ is said to be \emph{quadratic holomorphic} iff
\begin{equation}\index{cross-ratio perserving maps}
  \label{eq:quadHolo}
  \frac{\(f(y)-f(x)\)\(f(y')-f(x')\)}{\(f(x)-f(y')\)\(f(x')-f(y)\)}=
  \frac{\(y-x\)\(y'-x'\)}{\(x-y'\)\(x'-y\)}.  
\end{equation}
A rhombic tiling gives rise to two sets of isoradial circle patterns:
a set of circles of common radius $\delta$, whose centers are the
vertices of $\Gamma$ and intersections are the vertices of $\Gamma^*$
and vice-versa. Two interesting families of cross-ratio preserving
maps are given by circle patterns with the same combinatorics and
intersection angles as one of these two circle patterns.
\index{circle patterns}
\begin{figure}[htbp]
\begin{center}\input{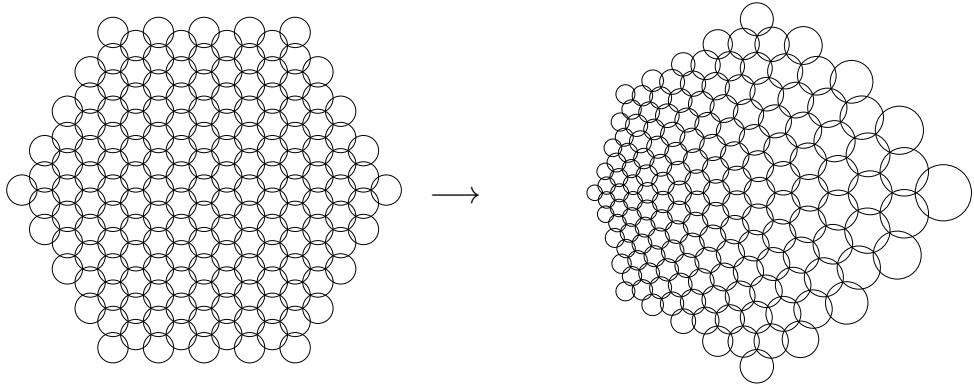}
\end{center}
\caption{A circle pattern with prescribed angles as a cross ratio preserving map.}  
\label{fig:trihexcirc}
\end{figure}
\begin{figure}[htbp]
\begin{center}\input{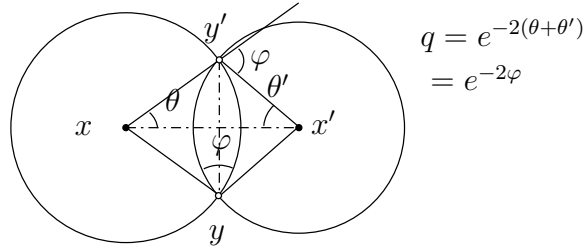}
\end{center}
\caption{The cross-ratio $q$ is given by the intersection angles.}  
\label{fig:quadAngles}
\end{figure}

A change of coordinates helps understanding diagonal ratio preserving
maps as a linearized version of the cross-ratio preserving maps. We
will say that the function $w:\diamondsuit_2\to\mathbb{C}$ solves the
\emph{Hirota system} if, around a face $(x,y,x',y')\in\diamondsuit_2$,
\begin{equation}\index{Hirota system}
  \label{eq:Hirota}
  (y-x)w(x)w(y)+(x-y')w(y')w(x)+(y'-x')w(x')w(y')+(x'-y)w(y)w(x') = 0.
\end{equation}
This is to be understood as a quadratic version of the Morera theorem
$\oint f\, dz=0$ and $w$ is a half of the derivative of
a holomorphic function:
\begin{propo}\label{th:Hirota}
  If $w$ solves the Hirota system, then the function
  $f:\diamondsuit_2\to\mathbb{C}$ defined up to an additive constant by
\begin{equation}
  \label{eq:fHirota}
  f(y)-f(x)=(y-x)w(x)w(y)
\end{equation}
is quadratic holomorphic.
\end{propo}
\begin{demo}{\ref{th:Hirota}}
  The function $f$ is well defined because the associated $1$-form is
  closed by definition of the Hirota system. The function $w$
  disappears in the cross-ratio of $f$, leaving the original cross-ratio.
\end{demo}
Conversely, a quadratic holomorphic function defines a solution to the
Hirota system, unique up to multiplication by $\lambda$ on $\Gamma$,
$1/\lambda$ on $\Gamma^*$. Concerning circle patterns families, $w$ is
real on the centers and unitary on the intersections, and encodes the
variation of radius, resp. of direction of the image of the circle:
\begin{equation}
  \label{eq:wwdZ}
  f(y)-f(x)=r(x)e^{i\theta(y)}(y-x).
\end{equation}

\begin{propo}
  \label{th:logDeriv}
  The logarithmic derivative of the Hirota system associated to a
  family of cross-ratio preserving maps is a diagonal ratio preserving
  map.
\end{propo}
In other words, for $(1+\epsilon g) w$ to continue solving the Hirota system at
first order, the deformation $g$ must satisfy
\begin{equation}
  \label{eq:epsg}
  \frac{g(y')-g(y)}{g(x')-g(x)}=
  \frac{f(y')-f(y)}{f(x')-f(x)}.
\end{equation}

\begin{demo}{\ref{th:logDeriv}}
The $\epsilon$ contribution of the closeness
condition~\eqref{eq:Hirota} for  $(1+\epsilon g) w$ gives
\begin{gather}
 (g(x)+g(y))\, f(x)\, f(y)\,  (y-x) + (g(y)+g(x'))\,f(y)\, f(x')\, (x'-y) \quad \quad \quad  \notag\\
\quad 
+ (g(x')+g(y'))\,w(x')\, w(y')\,  (y'-x') +
  (g(y')+g(x))\,w(y')\, w(x)\, (x-y') = 0,  \notag\\[.5cm]
\text{which reads, refering to $f$: }
  \qquad\frac{g(y')-g(y)}{g(x')-g(x)} = \frac{f(y')-f(y)}{f(x')-f(x)}.\notag
\end{gather}
\end{demo}

\subsection{Baecklund transformation}
The Baecklund transformation is a way to associate, to a given
solution of an integrable problem, a family of deformed solutions. The
two problems under consideration here are the linear and quadratic
holomorphicity constraints on each face. They are given by a linear,
resp. quadratic algebraic relation between the four values of a
solution at the vertices of each face. These relations involve only
values supported by the edges of the rhombus, which are equal on
opposite sides, namely the complex label $y-x=x'-y'$.

The Baecklund transformation is defined by imposing such constraints
over new virtual faces added over each edge, with ``vertical edges''
labelled by a complex constant $\lambda$:
\begin{defi}\index{Baecklund transformation}
  Given a linear holomorphic function $f\in\Omega(\diamondsuit)$,
  complex numbers $u,\lambda\in\mathbb{C}$, its Baecklund
  transformation  $f_\lambda=B_\lambda^u(f)$ is defined by
  \begin{eqnarray}
    \label{eq:BaeckDefi}
    f_\lambda(0)&=&u,\notag\\
\frac{f_\lambda(x)-f(y)}{f_\lambda(y)-f(x)}&=&\frac{\lambda+x-y}{\lambda+y-x}.
  \end{eqnarray}
  Given a quadratic holomorphic function $f$,
  complex numbers $u,\lambda\in\mathbb{C}$, its Baecklund
  transformation  $f_\lambda=B_\lambda^u(f)$ is defined by
  \begin{eqnarray}
    \label{eq:BaeckDefiQuad}
    f_\lambda(0)&=&u,\notag\\
\frac{f_\lambda(y)-f_\lambda(x)}{f_\lambda(x)-f(x)}
\frac{f(x)-f(y)}{f(y)-f_\lambda(y)}&=&
\frac{(y-x)^2}{\lambda^2}.
  \end{eqnarray}
\end{defi}
The right hand sides are the values respectively of the diagonal ratio
and cross-ratio of a parallelogram faces of sides $(y-x)$ and
$\lambda$ seen as ``over'' the edge $(x,y)\in\diamondsuit_1$.
\begin{figure}[htbp]
\begin{center}\input{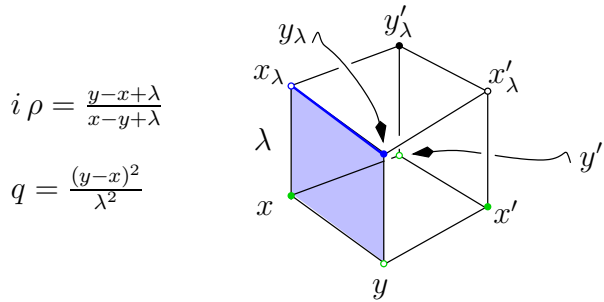}
\end{center}
\caption{The face $(x_\lambda,y_\lambda,x'_\lambda,y'_\lambda)$
  ``over'' the face $(x,y,x',y')\in\diamondsuit_2$.}  \label{fig:backlund}
\end{figure}

\begin{propo}
  This transformation is well defined in the critical case.
\end{propo}
This condition, called \emph{three dimensional consistency} is an
overdetermination constraint: if the cube ``over'' the face
$(x,y,x',y')\in\diamondsuit_2$ is split into two hexagons along the
cycle $(y,x',y',y'_\lambda, x_\lambda,y_\lambda)$, one can see that,
given values at these six vertices, the values at the centers of each
hexagons, namely at $x'_\lambda$ and at $x$ are overdetermined.
\begin{figure}[htbp]\index{three dimensional consistency}
\begin{center}\input{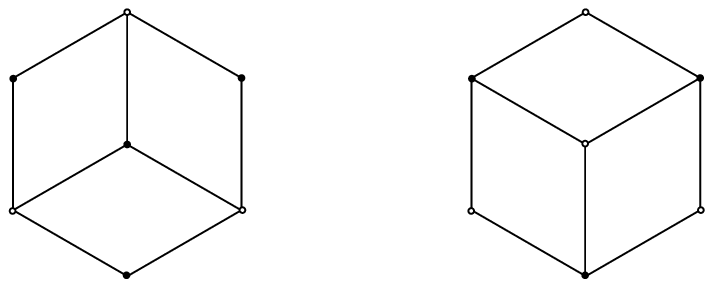}
\end{center}
\caption{The cube split into two hexagons yielding equivalent
  compatibility constraints.}  \label{fig:YB}
\end{figure}

Therefore only certain values at the six vertices are allowed, defined
by two algebraic relations between them. The compatibility condition
is that these two algebraic relations are equivalent. It is a simple
computation to check it is the case for critical maps.

This transformation verifies
\begin{equation}
  \label{eq:Baeck}
  B_{\lambda^{-1}}^{f(O)}\(B_\lambda^u(f)\) = f
\end{equation}
for any $(u,\lambda)$.  It is an analytic transformation in all the
parameters therefore its derivative is a linear map between the
tangent spaces, that is to say between diagonal ratio preserving maps,
\begin{equation}
  \label{eq:dBlambda}
   d\, B_\lambda^u(f):\Omega(f)\;\to\;\Omega\(B_\lambda^u(f)\).
\end{equation}
It is not injective and I define the discrete exponential at $f$ as being the
direction of this $1$-dimensional kernel. It can be characterized as a
derivative with respect to the initial value at the origin:
\begin{equation}
  \label{eq:defExpF}
  \exp_u({:}\lambda{:} f) := \frac{\partial\,}{\partial v}
  B_{\lambda^{-1}}^{v}\(B_\lambda^u(f)\)|_{v = f(O)}\; \in\; 
\ker\(d\, B_\lambda^u(f)\)
\end{equation}
because $B_\lambda^u\(B_{\lambda^{-1}}^{v}(g)\) = g$ for all $\lambda, g$
and $v$.\begin{figure}[htbp]
\begin{center}\input{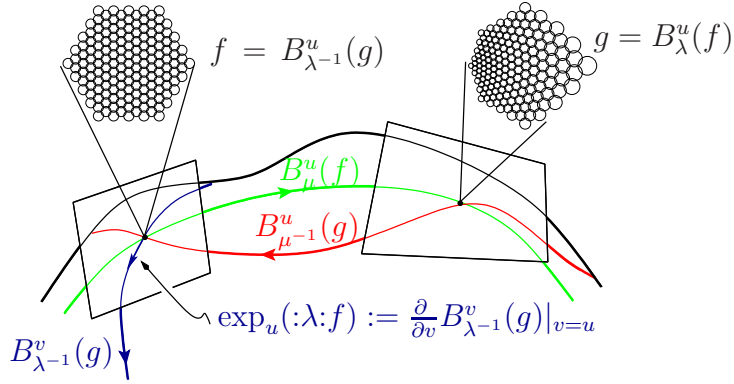}
\end{center}
\caption{ The discrete exponential $\exp_u({:}\lambda{:} f)$
  is the kernel of the linear transformation $d\, B_\lambda^u(f)$ (here
  $u=f(O)$).}         
\label{fig:tanBacklund}
\end{figure}

As in the discrete exponential case, the value of the Baecklund
transformation at a given vertex is the image of values at
neighbouring vertices by a homography. These homographies can be
encoded as projective operators $L(e;\lambda)\in GL_2(\mathbb{C})[\lambda]$
lying on the edges $e\in\diamondsuit_2$, called a \emph{zero curvature
representation}:
\begin{equation}\label{L CR}
\renewcommand{\arraystretch}{1.4}
L((x,y);\lambda)=\begin{pmatrix}
\lambda+y-x & -2(y-x)(f(x)+f(y)) \\ 0 & \lambda+x-y \end{pmatrix}
\quad \text{ for the linear case,}
\end{equation}

\begin{equation}\label{L Hirota}
\renewcommand{\arraystretch}{1.4}
L((x,y);\lambda)=\begin{pmatrix} 1 & -(y-x) w(y) \\
-\lambda(y-x)/w(x) & w(y)/w(x) \end{pmatrix} \quad \text{ for the
Hirota system.}
\end{equation}
Then we define~\cite{NRGO} the \emph{moving frame} $\Psi:\diamondsuit_2\to
GL_2(\mathbb{C})(\lambda)$ by a prescribed value at the origin and
recursively by $\Psi(y;\lambda) = L((x,y);\lambda) \Psi(x; \lambda)$
and its logarithmic derivative with respect to $\lambda$\index{moving frame}
\begin{equation}\label{A}
  A(e;\lambda)=\frac{d\Psi(e;\lambda)}{d\lambda}\Psi^{-1}(e;\lambda)
\end{equation}
is meromorphic in $\lambda$ for each edge $e$. We call $f$, resp. $w$
\emph{isomonodromic} if the positions and orders of the poles don't
depend on the edge $e$. The two points discrete Green function (the discrete
logarithm) $G(O,x)$, inverse of the Laplacian in the sense that
\begin{equation}\index{discrete Green function}
  \label{eq:Green}
  \Delta G(O,\bullet) = \delta_{O,\bullet}
\end{equation}
can be constructed as the unique isomonodromic solution with some
prescribed data~\cite{M0402097}, which allows us to give an explicit
formula for it, recovering results of Kenyon~\cite{Ken02I}: an
integral over a loop in the space of discrete exponentials:
\begin{equation}
  \label{eq:KenyonGreen}
  G(O,x)=-\frac{1}{8\,\pi^2\,i}\oint_{C}\,
\exp({:}\lambda{:}\,x) \, \frac{\log \frac{\delta}{2} \lambda}{\lambda}\, d\lambda
\end{equation}
where the integration contour $C$ contains all the possible poles of
the rational fraction $\exp({:}\lambda{:}\,x)$ but avoids the half
line through $-x$. It is real (negative) on half of the vertices and
imaginary on the others.  Because of the logarithm, this imaginary
part is multivalued.

\printindex

\bibliographystyle{unsrt} \bibliography{these}
\end{document}